\documentclass{article}
\usepackage[latin1]{inputenc}
\usepackage[english]{babel}
\usepackage{amsthm,amssymb,amsmath,epsfig,graphics}
\usepackage{color}
\usepackage[all]{xy}
\usepackage{mathrsfs,bbm}
\usepackage{graphics}
\usepackage{accents}
\usepackage{titlesec}
\titleformat{\subsection}[runin] {\normalfont\it} {\thesubsection}{.5em}{}
\usepackage[colorlinks,
            linkcolor=blue,
            anchorcolor=blue,
            citecolor=blue
            ]{hyperref}
\newfam\msbmfam\font\tenmsbm=msbm10\textfont
\msbmfam=\tenmsbm\font\sevenmsbm=msbm7
\scriptfont\msbmfam=\sevenmsbm
\def\NN{\mathbb N}
\def\ZZ{\mathbb Z}

\def \RR{\mathbb R}

\def\EE{{\bb E}}
\def\EE{{\mathbb E}}
\def\QQ{\mathbb Q }

\def\proof{{\noindent \it Proof\quad  }}

\def\rme{{\rm
e}}

\def\be{{\beta}}
\def\de{{\delta}}
\def\la{{\lambda}}
\def\si{{\sigma}}

\def\rZ{{\cal Z}}
\def\De{{\Delta}}

\def\Om{{\Omega}}\def\al{{\alpha}}
\def\be{{\beta}}\def\Ga{{\Gamma}}
\def\ga{{\gamma}}\def\de{{\delta}}\def\De{{\Delta}}

\def\si{{\sigma}}

\def\la{{\lambda}}
\def\vare{{\varepsilon}}
\def\vr{{\varrho}}

\def\i{{\imath}}

\newtheorem{thm}{Theorem}[section]
\newtheorem{lemma}[thm]{Lemma}

\newtheorem{prop}[thm]{Proposition}

\newtheorem{remark}[thm]{Remark}

\def\ZZ{{\mathbb Z}}
\def\PP{{\mathbb P}}

\def\al{\alpha}

\def\Ad2{{| A-\widetilde A |_{L^{2}_{loc}} }}

\def\om{{\omega}}\def\Om{{\Omega}}

\def\rZ{{\cal Z}}
 
\def \eref#1{\hbox{(\ref{#1})}}

\newcommand{\D}{{\mathrm d}}

\newcommand{\nc}{\newcommand}
\nc{\bi}{\bibitem}
\newcommand{\e}{{\rm e}}

\begin{document}
\title{From Directed Polymers in Spatial-correlated Environment to Stochastic Heat Equations Driven by Fractional Noise in $1+1$ Dimensions}
\author{
Guanglin  {\sc  Rang}\\
{\it\footnotesize School of Mathematics and Statistics, Wuhan University, Wuhan 430072,China}\\
{\it\footnotesize Computational Science Hubei Key Laboratory, Wuhan University, Wuhan, 430072, China}\\
{\it\footnotesize glrang.math@whu.edu.cn}}
\date{}
\maketitle
\begin{abstract}  We consider the limit behavior of partition function of directed polymers in random environment, which is represented by a linear model instead of a family of i.i.d.variables in $1+1$ dimensions. Under the assumption on the environment that its spatial correlation decays algebraically, using the method developed in [Ann. Probab., 42(3):1212-1256, 2014], we show that the scaled partition function, as a process defined on $[0,1]\times\RR$,  converges weakly to the solution to some stochastic heat equations driven by fractional Brownian field. The fractional Hurst parameter is determined by the correlation exponent of the random environment. Here multiple It\^{o} integral with respect to fractional Gaussian field and spectral representation of stationary process are heavily involved.\newline
 {\bf Key words}: stochastic heat equation, random walk, partition function, fractional noise, stationary fields, multiple It\^{o} integral  \newline
{\bf 2010 MR Subject Classification}:  60F05,60H15,82C05.
\end{abstract}
\section{Introduction}

We are concerned with the large time behavior of directed polymers in random environment with long-range correlation.  The directed polymers is used to model the behavior of a polymer chain when it stretches in some media with impurities or charges. It consists of   a directed random walk and a family of random variables attached in every space-time point representing the random environment. Mathematically, let ${S}=\{S_k\}_{k=0}^\infty$ be a nearest-neighbor path starting from the origin in $\ZZ^d$,  and $\om=\{\om(i,x), (i,x)\in \NN\times\ZZ^d\}$ be a family of real-valued random variables appearing as the environment. In what follows $\PP$ and $\QQ$ are the probability measures corresponding to the random walk $S$ and the environment variables $\om$, and the notations ${\mathbb E}_{\mathbb P},{\mathbb E}_{\mathbb Q}$ represent the expectation with respect to ${\mathbb P},{\mathbb Q}$, respectively. 

Given a fixed environment $\om$, the $n$-step energy of a path $S$ is 
$$
H^\om_n(S)=\sum_{i=1}^n\om(i,S_i),
$$
and its random polymer measure is given as the usual Gibbsian form by 
$$
\PP^\om_n(S)=\frac{1}{Z_n^\om(\be)}\rme^{\be H_n^\om(S)}\PP(S),
$$
where $\be$ is the inverse temperature and 
\begin{align}\label{partition1}
Z_n^\om(\be)=\sum_S{\rm e}^{\be H_n^\om(S)}\PP(S)=\EE_\PP({\rm e}^{\be H_n^\om(S)})
\end{align}
is the partition function. Also, for any $x\in\ZZ$, we have the $0$-to-$x$ partition function defined by
\begin{align}\label{partition2}
Z_n^\om(x;\be)=\sum_S{\rm e}^{\be H_n^\om(S)}\PP(S|S_n=x)=\EE_\PP({\rm e}^{\be H_n^\om(S)}\mathbbm{1}_{\{S_n=x\}}).
\end{align}

Because of the interaction between the chain and the environment a basic phenomenon of the directed polymer is that the average path wanders significantly farther along the axes transverse to the directed axis (time direction) than the purely entropic spreading walk. In other words, the transverse fluctuations in random media are supperdiffusive, i.e.,  $\overline{\langle x^2(t)\rangle}\sim t^{2v}$ with $v\geq 1/2$, differing from the standard diffusion behavior in no-random media, here the sharp brackets represent the average over all the paths in configuration space while the bar means the average under the environment. The constant $v$ is usually referred as diffusive exponent, to which is connected other exponents, such as dynamic exponent $z=\frac{1}{v}$ and fluctuations exponent $\chi$ with $\chi=2v-1$. These relationships, among others, have been established analytically (replica method, renormalization group analysis) or observed numerically in some models in the setting of the environment being a collection of i.i.d. random variables or in continuum contexts. For more comments see \cite{alberts2014,PhysRevLett.56.889,Lacoin2011,Medina1989Burgers} and references therein.  

The positive parameter $\beta$ and the dimension $d$  play  important role in determining  the behavior of the polymer. $\be=0$ is the usual random walk, while $\be=\infty$ corresponds to the last passage percolation(\cite{Johansson2000Shape}).  There is a critical value $\beta_c$ such that in the phase $(0,\beta_c)$, called weak disorder, entropy dominates, which means that the trajectory of the polymer conserves all the essential
features of the nondisordered model (i.e., simple random walk); whereas in the phase $(\beta_c,\infty)$, called strong disorder, energy dominates, which means that the end position of a trajectory up to time $t$ to the origin is greater than $t^{\frac{1}{2}}$ or in other words, the chain tends to go farther from the origin to reach a more favorable environment. For more reviews about this subject see the introductory survey of the Kardar-Parisi-Zhang equation (KPZ) given by Quastel (\cite{Quastel2011KPZ}). In \cite{Alberts-Clark2017, alberts2014}, some new scalings of $\beta$ with the length of the polymer were proposed (the object in \cite{alberts2014} is as the same as this work, while in \cite{Alberts-Clark2017} is for a directed polymer model defined on a hierarchical diamond lattice). These new scalings lead to an intermediate disorder regime, which sits between weak and strong disorder regime. Speaking roughly, under this scaling and diffusive scaling as well, the conditioned polymer measure converges in law to a random field, to which a solution to stochastic heat equation with multiplicative time space white noise is related. Specifically, under the assumption that $\{\om(i,x), (i,x)\in \NN\times\ZZ^d\}$ are independent identical distribution (IID for short) with the existence of appropriate moments, the modified partition function (see \eqref{mpartition1} below) of the original exponential form \eqref{partition1} was expanded as an $n$-order polynomial of $\beta$, then every monomial was proved to converge weakly to a multiple Wiener-It\^{o} integral of transition probability of simple symmetrical random walk with respect to time space white noise. Hence the scaled modified partition function and moreover the original one  converge weakly to the solution to the equation
\begin{align}\label {SHE}
\frac{\partial u(t,x)}{\partial t}=\frac{1}{2}\Delta u(t,x)+\sqrt{2}\beta u(t,x){\eta}(t,x)
\end{align}
since it admits a chaos solution in terms of multiple It\^{o} integral, where $\eta(t,x)$ is time space white noise, i.e.,
\begin{align*}
\overline{\eta(t,x)\eta(t',x')} =\delta(x-x')\delta(t-t')
\end{align*}
and $\delta(x)$ is Dirac function. For more detail about multiple It\^{o} integral and stochastic heat equation (SHE for short) see \cite{Hu2001Heat,Hu2002Chaos,Hu2015Stochastic} or Section 3 below, for general stochastic heat equation (Anderson model) see \cite{SONG201737}. This builds a bridge between discrete polymer models and stochastic heat equations or other stochastic growth models, e.g., KPZ equations, whose fluctuation characteristics are shown to be more universal.  This idea also inspired a recent work \cite{CRN2017}, in which a united framework is developed for studying the continuum and weak disorder scaling limits of statistical mechanics systems that are disorder relevant, including the disordered pinning model (\cite{CRN2016}), the stable directed polymer model in dimension 1 + 1, and the two-dimensional random field Ising model. 

In this paper we are focused on the limit behavior of random polymer in a spatial correlated random environment in $1+1$ dimension. Concretely,  we assume that the environment $\om=\{\om(n,x), (n,x)\in \NN\times\ZZ\}$ has the following form:
\begin{align}\label{envir1}
\om(n,x)=\sum_{y=-\infty}^\infty\psi_{y-x}\xi_{n,y},
\end{align}
with $\psi_j\sim\de|j|^{-\al}$ and $1/2<\al<1, \delta>0$. Here $\{\xi_{i,j}:i\in \NN, j\in \ZZ\}$ be a family of independent identical distribution variables with $\EE_{\mathbb Q} \xi_{i,j}=0$ and $\EE_{\mathbb Q} \xi^2_{i,j}=1$ for any $i,j$; $\QQ$ is the probability measure corresponding the environment variables. This is a good setting for us since it can be dealt with  to some extent and does not lose generality. 

Since $1/2<\al<1$, the stationary field $\om$ has long-memory property, hence  it tends to some Gaussian fractional field $\dot{W}(t,x)$ (see Section 3 below) with Hurst parameter $1/2<H<1$ after a time space scaling. Let $\Lambda (\cdot)$ be the Log-Laplace of $\om$ and $\varrho=\frac{H}{2}$. Then by a similar procedure as in \cite{alberts2014} we can show the scaled partition function 
\begin{equation}\label{C11}
\frac{\sqrt{n}}{2}\e^{-nt\Lambda(\be n^{-\vr})}Z^\om_n(nt,\sqrt{n}x;\be n^{-\vr})\Longrightarrow u(t,x) \quad \text{weakly}
\end{equation}
in the sense of process level, where $u(t,x), (t,x)\in [0,1]\times\RR$, is the solution to  
the following equation  
\begin{align}\label {SHE1}
\frac{\partial u(t,x)}{\partial t}=\frac{1}{2}\Delta u(t,x)+\sqrt{2}\beta u(t,x)\dot{W}(t,x)
\end{align}
with initial data $u(0,x)=\delta(x)$. See Theorem \ref{convergence12} in Section 4.  

For this sake, we first consider the convergence of the modified partition function $\mathfrak {Z}_n^\om(\be)$ obtained by using $1+x$ to take place of $\e^x$ in the partition function $Z_n^\om(\be)$, i.e.,
\begin{align}
\label{mpartition1}
\mathfrak {Z}_n^\om(\be)=&{\mathbb E_{\mathbb P}}\left[\prod_{i=1}^n(1+\be\om(i,S_i))\right]\nonumber\\
=&{\mathbb E_{\mathbb P}}\left[1+\sum_{k=1}^n\be^k\sum_{{\mathbf i}\in D_k^n}\prod_{j=1}^k\om(i_j,S_{i_j})\right].
\end{align}
See \eqref{porder1} for the notation $D_k^n$ . 

Formally, we can see that every term in the above equation is the discrete multiple stochastic integral of the transition probability of the random walk $S$ with respect to $\om$. 
Then, after a $\be$ scaling, we show that every term in the following equation
\begin{align}
\label{mpartition11}
\mathfrak {Z}_n^\om(\be n^{-\vr})
=&{\mathbb E_{\mathbb P}}\left[1+\sum_{k=1}^n\be^kn^{-k\vr}\sum_{{\mathbf i}\in D_k^n}\prod_{j=1}^k\om(i_j,S_{i_j})\right].
\end{align}
converges to the corresponding multiple stochastic integral of the transition density of Brownian motion with respect to the fractional Gaussian field $\dot{W}(t,x)$. That is, 
the $k$-order functionals \eqref{Uk}, called U-statistics, of the transition density of symmetric random walk,  converges to $k$th-multiple It\^{o} integral with respect to the fractional field.
We know that the multiple It\^{o} integral with respect to fractional Brownian motion is defined by Hermite function or Wick product of Wiener integrals instead of the linear extension of integral of indicators of $k$-dimension rectangle as in Gaussian white noise case. This leads to considerable computation in this work. We first show that it is true for $k=1$, then transit to the case of $k\geq2$ by resorting to the recursive identities \eqref{contract1} and \eqref{contract2} for multiple It\^{o} integral. This is contained in Theorem \ref{mainthm}, which says that this convergence holds for all symmetry function $f\in {\mathcal L}_H^{\otimes k}$ (see \eqref{ksym}). On the base of Theorem \ref{mainthm} we have the intermediate result Theorem \ref{convergence}, which claims that the scaled point-to-point modified partition function
$$
\frac{\sqrt{n}}{2}\mathfrak {Z}_n^\om(nt,\sqrt{n}x;\be n^{-\vr})\longrightarrow u(t,x) \quad \mbox{weakly, as}\quad n\longrightarrow\infty,
$$
where $u(t,x)$ is the solution to equation \eqref{SHE1}.

In order to go back to the weak convergence of the true partition function from that of the modified partition function, we define a non-linear functional 
\begin{equation}\label{omegan}
{\om}_n(i,x)=\frac{\e^{\be n^{-\vr}\om(i,x)-\Lambda(\be n^{-\vr})}-1}{\be n^{-\vr}}\stackrel{\De}{=}F(n,\om(i,x))
\end{equation}
of $\om$. Since $F(n,\om(i,x))$ in essence is $\om$ when $n$ is big enough, we can use 
the theory about the central limit theorem for non-linear functionals of stationary process developed in \cite{Ho and Hsing1998Annals} to show that the corresponding scaled U-statistics e.g., 
\begin{align*}
\mathfrak {Z}_n^{{\om}_n}(\be n^{-\vr})=&1+\sum_{k=1}^n\be^kn^{-k\vr}\sum_{{\mathbf i}\in D_k^n}\sum_{{\mathbf x}\in {\ZZ}^k}\left[\prod_{j=1}^k{\om}_n(i_j,x_j)p_k({\mathbf i},{\mathbf x})\right]
\end{align*}
has a weak limit (see Theorem \ref{omeganconvergence}). Since
\begin{align*}
\e^{-n\Lambda(\be n^{-\vr})}{\mathbb Z}^\om_n&={\EE}_{\PP}\Pi_{i=1}^n(1+\be n^{-\vr}{\om}_n(i,S_i))\\
&=1+\sum_{k=1}^n\be^kn^{-k\vr}\sum_{{\mathbf i}\in D_k^n}\sum_{{\mathbf x}\in {\ZZ}^k}\left[\prod_{j=1}^k{\om}_n(i_j,x_j)p_k({\mathbf i},{\mathbf x})\right],
\end{align*}
all things above can help us arrive at the main result, Theorem \ref{convergence12}  or see \eqref{C11}. By virtue of the spatial correlation of the environment, spectral representation of the correlation is used extensively to help us estimate the variance of U-statistics (see \eqref{clt1} below) and prove the tightness of approximation process, Theorem \ref{tight1}, Section 5.

Before closing this section we remark some facts about the continuum random directed polymers with environment being correlated in space and/or in time.  Although the correlated environment in discrete setup is considered first in this paper, the study for the continuum case can be traced in \cite{Medina1989Burgers}. There the time space correlation of noise $\dot{W}(t,x)$ is formulated as 
\begin{align*}
\overline{\dot{W}(t,x)\dot{W}(t',x')}\sim |x-x'|^{2\rho-d'}|t-t'|^{2\theta-1}
\end{align*}
with $d'=d-1$, different from the one in equation \eqref{SHE1}. Then the solutions $h(t,x)$, called height functions, to the following KPZ equations driven by random forcing $\dot{W}(t,x)$
\begin{align}\label {KPZ1}
\frac{\partial h(t,x)}{\partial t}=\nu\nabla^2h(t,x)+\frac{\lambda}{2}(\nabla h(t,x))^2+\dot{W}(t,x)
\end{align}
satisfy
\begin{align*}
\overline{|h(t,x)-h(t',x')|^2}\sim|x-x'|^{2\zeta}f\left[\frac{|t-t'|}{|x-x'|^z}\right]\end{align*} 
with some functions $f$, where $\zeta$ is roughening exponent, $z$ is dynamic exponent. (Throughout this paper, the symbol $\sim$ means that the ratio of the quantities lying its two sides goes to 1 as the corresponding argument tending to some limit. Here $\sim$ means that the arguments $|x-x'|$ and $|t-t'|$ are large.)  Renormalization analysis and numerical simulation indicate that in $d'=1$, $\zeta, z$ vary with $\rho\in[0,1]$,  and large $\rho$ tends to roughen the interface or increase the roughening exponent in the absence of temporal correlation (more complicated when temporal correlation appears). 
Recently in \cite{Lacoin2011}, a Brownian directed polymer with spatial correlation was studied. In that model, Brownian motion $B$ and a real centered Gaussian field $\om$ were in place of the random walk $S$ and IID environment in \eqref{partition1}, respectively. Its Hamiltonian and Polymer measure are given by
$$
H _ { \omega , t } ( B ) = H _ { t } ( B ) : = \int _ { 0 } ^ { t } \omega \left( \mathrm { d } s , B _ { s } \right)
$$
and
$$
\mathrm { d } \mu _ { t } ^ { \beta , \omega } ( B ) : = \frac { 1 } { Z _ { t } ^ { \beta , \omega } } \exp \left( \beta H _ { t } ( B ) \right) \mathrm { d } P ( B )
$$
with partition function $Z _ { t } ^ { \beta , \omega } : = P \left[ \exp \left( \beta H _ { t } \right) \right]$.
The strong spatial correlation with polynomial decay is assumed in that context, by the technique from stochastic analysis, the behavior of polymer in terms of free energy, fluctuation exponent and volume exponent are considered systematically. The similarity among these models stimulates us to consider in this work the discrete approximation, by which we wish find some important things, e.g., the universality constants different from that usual KPZ equations have shared. 

Here is the structure of this paper. In Section 2, some notations of the transition function of random walk is given as well as a basic central limit theorem, whose proof is a little 
different from that in Section 4. In Section 3 some basic facts about fractional Gaussian fields and stochastic heat equations
with multiplicative fractional noise are listed. In Section 4 we state and prove the main results Theorem \ref{mainthm} in this paper; we take a large space on the proof that the scaled partition function converges  to the solution to some stochastic heat equation driven by fractional Gaussian fields in the sense of finite dimensional distribution (Theorem \ref{convergence} and Theorem \ref{convergence12}). Finally, in Section 5 we prove that the approximation process is tight via checking Kolmogorov's criterion.  

\section{Some notations and a central limit theorem}

For the free random walk $S$ on $\ZZ$, we denote by ${\mathfrak S}_n$ its $n$-step configuration space, a subset  consisting of all possible discrete  paths started at zero, i.e., ${\mathfrak S}_n=\{(0,x_1,\dots,x_n):x_i\in{\mathbb Z}, |x_i-x_{i+1}|=1, i=0,1,\dots,n-1,x_0=0\}$. Let $p(n,x)={\mathbb P}(S_n=x)=2^{-n}\binom{n}{(n+x)/2}$ if $n\leftrightarrow x$, otherwise 0.
Here $n\leftrightarrow x$ means that $n$ and $x$ have the same parity. Let $[n]=\{1,2,\dots,n\}$ and, for $k\in[n]$,
\begin{align}\label {porder1}
D_k^n=\{{\textbf i}=(i_1,i_2,\dots,i_k)\in [n]^k: 1\leq i_1<i_2<\dots<i_k\leq n\}.
\end{align}
For ${\textbf i}\in D_k^n$, write $S_{\textbf i}=(S_{i_1},S_{i_2},\dots,S_{i_k})$.
Then we have the $k$-dimension joint probability
\begin{align}\label{k-transition}
p_k({\mathbf i},{\mathbf x}):={\mathbb P}(S_{\textbf i}={\textbf x})={\mathbb P}(S_{i_1}=x_1,S_{i_2}=x_2,\dots,S_{i_k}=x_k)
=\prod_{j=1}^kp(i_j-i_{j-1},x_j-x_{j-1})
\end{align}
for ${\textbf x}=(x_1,x_2,\dots,x_k)\in {\mathbb Z}^k$ and $i_0=0,x_0=0$. Here, it necessitates that $i_j\leftrightarrow x_j$ for all $j=1,2,\dots,k$, which is denoted by ${\mathbf i}\leftrightarrow{\mathbf x}$. 
 
 We can also extend $p_k({\mathbf i}, \cdot)$ to all of ${\mathbb R}^k$ by defining a density function by (\cite{alberts2014})
 $$
 \bar{p}_k({\mathbf i}, {\mathbf x})=2^{-k}p_k({\mathbf i}, [{\mathbf x}]_{\mathbf i}),
 $$
where the jth component of $[{\mathbf x}]_{\mathbf i}\in{\mathbb Z}^k$, denoted by $([{\mathbf x}]_{\mathbf i})_j=[x_j]_{i_j}$, is the closest integer to $x_j$ in ${\mathbb Z}$ such that $i_j\leftrightarrow [x_j]_{i_j}$.  It is known that $\bar{p}_k({\mathbf i}, {\mathbf x})$ is the finite dimensional distribution for random walk $X_n=S_n+U_n$, where $\{U_n,n\in{\mathbb N}\}$ is a family of uniform IID random variables in $(-1,1)$. 

Furthermore, for $k,n\in{\mathbb N}$, define ${p}^n_k$ by 
\begin{equation}\label{SP1}
{p}^n_k({\mathbf t},{\mathbf x})=\bar{p}_k(\lfloor n{\mathbf t}\rfloor, \sqrt{n}{\mathbf x}){\mathbf1}_{\lfloor n{\mathbf t}\rfloor\in D_k^n},
\end{equation}
for $({\mathbf t}, {\mathbf x})\in [0,1]^k\times{\mathbb R}^k$, where $\lfloor n{\mathbf t}\rfloor=(\lfloor nt_1\rfloor,\lfloor nt_2\rfloor,\cdots,\lfloor nt_k\rfloor)$ and $\lfloor nt\rfloor$ is the largest integer not exceeding $nt$.

Let $\om=\{\om(i,x): i\in\{0\}\cup\NN,x\in\ZZ\}$ be the form of \eqref{envir1}. Then, one has (see \cite{Hosking1996261})
$$
\EE(\om(i,x)\om(j,y))=\de_{ij}\ga(x-y),
$$
where $\de_{ij}$ is Kronecker and $\ga(k)\sim \la |k|^{1-2\al}$ for large integer $k$ and $\la=\de^2\frac{\Ga(2\al-1)\Ga(1-\al)}{\Ga(\al)}$. 

For later use, let $G(\D \eta)$ be the spectral measure of the correlation function $\ga$, i.e.,
\begin{align}\label{spectrum}
\ga(k)=\int_{-\pi}^\pi\e^{\i k\eta}G(\D\eta), ~~~\forall k\in \ZZ,
\end{align}
with $\i^2=-1$.
\begin{remark}\label{coefficient}
{\rm Taking $H=3/2-\al$, then $1/2<H<1$, we call it Hurst index, with which a fractional Gaussian field will be involved. In what follows we may transit between the two parameters freely unless otherwise stated. We can also adjust the coefficient $\de$ such that $\la=H(2H-1)$. } 
\end{remark}

For every $N\in\NN$, we define a new  measure $G_N$ by
$$
G_N(A)=N^{\al-1/2}G(N^{-1/2}A),~~~~~~~A\in \mathcal{B}(\RR).
$$
Then, from \cite[Proposition 1]{Dobrushin1979Non}, there exists a locally finite measure $G_0$, such that 
$$
\lim_{N\rightarrow\infty}G_N=G_0
$$
in the sense of locally weak convergence. Furthermore, $G_0$ has a spectral density $D^{-1}|\eta|^{1-2H}$ with $D=2\Ga(2-2H)\cos(1-H)\pi$. Therefore  $G_0$ has a density as follows:
\begin{equation}\label{limitofspectrum}
G_0(\D x)=\frac{D^{-1}}{2\pi}\int_{\RR}\e^{\i x\eta}|\eta|^{1-2H}\D \eta\D x=\frac{\D x}{2\pi|x|^{2-2H}}.
\end{equation}
Let
\begin{equation}\label{K}
K(x)=\frac{H(2H-1)}{|x|^{2-2H}},
\end{equation}
which is the covariance of the increment of some fractional Brownian motion with Hurst parameter $H>1/2$ (\cite{reed1995spectral}). Then we have $G_0(\D x)=\frac{K(x)\D x}{2H(2H-1)\pi}$.

We are going to scale the partition functions and consider their limit behavior, which is postponed to Section 4.  In this section we show in the correlated setting a CLT holding for weighted sum, which is of interest on its own right. The proof will be mentioned repeatedly later, say, \eqref{lindeberg2}.

\begin{prop}\label{clt1}
Let $\om$ be given by \eqref{envir1}, $S$ be the symmetrical random walk on ${\mathbb Z}$ started at the origin, and let  $\varrho=H/2$. Then
\begin{align}\label{CLT1}
\beta n^{-\varrho}\sum_{i=1}^n\sum_{x\in\ZZ}\om(i,x)\PP(S_i=x)\stackrel{\cal D}{\longrightarrow}N(0,\si^2)
\end{align}
with $\si^2=\frac{\be^2\Ga(H-1/2)}{(3-2H)\pi}$.
\end{prop}
\proof 
First, we compute the variance of the expectation of the $n$-step energy with respect to probability measure ${\mathbb Q}$
\begin{align*}
A_n^2:=&\EE_{\QQ}(\EE_{\PP}( H^\om_n(S)))^2
=\EE_{\QQ}\bigg(\sum_{i=1}^n\sum_{x\in\ZZ}\om(i,x)\PP(S_i=x)\bigg)^2\\
=&\sum_{i=1}^n\EE_{\QQ}\bigg[\sum_{x\in\ZZ}\om(i,x)\PP(S_i=x)\bigg]^2
=\sum_{i=1}^n\EE_{\QQ}\bigg[\sum_{k=0}^i\om(i,2k-i)\PP(S_i=2k-i)\bigg]^2\\
=&\sum_{i=1}^n\sum_{k,l=0}^i\EE_{\QQ}\bigg[\om(i,2k-i)\om(i,2l-i)\PP(S_i=2k-i)\PP(S_i=2l-i)\bigg]\\
=&\sum_{i=1}^n\sum_{k,l=0}^i2^{-2i}\ga(2k-2l)\binom{i}{k}\binom{i}{l}.
\end{align*}
By \eref{spectrum}, one knows that the above quantity equals
\begin{align}\label{E1}
\begin{split}
&\sum_{i=1}^n\int_{-\pi}^{\pi}\sum_{k,l=0}^i2^{-2i}\binom{i}{k}\binom{i}{l}\e^{2\i k\eta}\e^{-2\i l\eta}G(\D\eta)\\
=&\sum_{i=1}^n\int_{-\pi}^{\pi}\bigg|\frac{\e^{\i\eta}+\e^{-\i\eta}}{2}\bigg|^{2i}G(\D\eta)=\int_{-\pi}^{\pi}\frac{|\cos\eta|^2-|\cos\eta|^{2n+2}}{1-|\cos\eta|^2}G(\D\eta)\\
=&n^{-1+H}\int_{-\sqrt{n}\pi}^{\sqrt{n}\pi}g_n(\eta)G_n(\D\eta),
\end{split}
\end{align}
where $g_n(\eta)=\frac{|\cos(\eta/\sqrt{n})|^2-|\cos(\eta/\sqrt{n})|^{2n+2}}{1-|\cos(\eta/\sqrt{n})|^2}\sim n\frac{1-\e^{-\eta^2}}{\eta^2}$. Put $g_0(\eta):=\frac{1-\e^{-\eta^2}}{\eta^2}$, which is integrable with respect to $G_0(\D\eta)$. Then we have, by \eref{limitofspectrum},
$$\si^2=\be^2\int_{\RR}g_0(\eta)G_0(\D\eta)=\be^2\int_{\RR} \frac{1-\e^{-\eta^2}}{2\pi|\eta|^{4-2H}}\D \eta=\frac{\be^2\Ga(H-1/2)}{(3-2H)\pi},\quad \text{
and}\quad 
A_n^2\sim n^H\si^2/\be^2.
$$

Now let 
\begin{align*}
&X_n=\sum_{i=1}^n\sum_{x\in{\cal E}_i}\om(i,x)a_{i,x}
=\sum_{y=-\infty}^\infty\sum_{i=1}^n\sum_{x\in{\cal E}_i}a_{i,x}\psi_{y-x}\xi_{i,y}
\end{align*}
with ${\cal E}_i=\{-i,-i+2,\cdots, i-2,i\}$ and non-random weight $a_{i,x}=\PP(S_i=x)$ satisfying $\sum_{x\in{\cal E}_i}a_{i,x}=1$ for every $i\in\NN$. Then 
$$
\EE_{\QQ}(X_n)^2=\sum_{y=-\infty}^\infty\sum_{i=1}^n\bigg(\sum_{x\in{\cal E}_i}a_{ix}\psi_{y-x}\bigg)^2=A_n^2\longrightarrow\infty
$$
as $n\rightarrow\infty$ by the previous calculation. 

We can use the method in \cite{Ibragimov1971Independent} or the corrected version \cite{Hosking1996261} to show $\frac{X_n}{A_n}\stackrel{\cal D}{\longrightarrow}N(0,1)$ . Putting 
$$\sum_{i=1}^n(b_{i,y})^2:=\sum_{i=1}^n\bigg(\sum_{x\in{\cal E}_i}a_{i,x}\psi_{y-x}\bigg)^2,
$$
then we have, for $i=1,2,\dots,n,$
\begin{align}\label{lindeberg1}
\begin{split}
b^2_{i,y}\leq&\sum_{x\in{\cal E}_i}a^2_{i,x}\sum_{x\in{\cal E}_i}\psi^2_{y-x}
\leq C\sum_{y=-\infty}^\infty\psi_y^2=C,
\end{split}
\end{align}
by Schwarz's inequality. Hence $b'_{n,i,y}:=\frac{b_{i,y}}{A_n}\leq\frac{C}{A_n}:=b_n\rightarrow 0$ uniformly in $y$ as $n\rightarrow\infty$. Thus 
$$
\frac{X_n}{A_n}=\sum_{y=-\infty}^\infty\sum_{i=1}^nb'_{n,i,y}\xi_{i,y}~~\text{with}~~
\sum_{y=-\infty}^\infty\sum_{i=1}^n(b'_{n,i,y})^2=1.
$$

Let $(\vare_n)$ be a sequence of positive number such that $\vare_n\rightarrow 0$ as $n\rightarrow \infty$. For each $n$, choose $N$ large enough so that 
$$
\sum_{|y|>N}\sum_{i=1}^n(b'_{n,i,y})^2<\vare_n,
$$
and consider the sequence of independent random variables $\{\eta_{n,0}\}\cup\{\eta_{n,i,y}; i=1,2,\dots,n, y=1,\dots 2N+1\}$, where
$$
\eta_{n,0}=\sum_{|y|>N}\sum_{i=1}^nb'_{n,i,y}\xi_{i,y}~~\text{and}~~\eta_{n,i,y}=b'_{n,i,-N+y-1}\xi_{i,-N+y-1}.
$$

Thus, one has, for any $\vare>0$
\begin{align*}
&\int_{|z|>\vare}z^2\D {\mathbb Q}(\eta_{n,0}<z)+\sum_{i=1}^n\sum_{y=1}^{2N+1}\int_{|z|>\vare}z^2\D {\mathbb Q}(\eta_{n,i,y}<z)\\
\leq &\vare_n+\sum_{i=1}^n\sum_{y=1}^{2N+1}\int_{|z|>{\vare}|b'_{n,i,-N+y-1}|^{-1}}(b'_{n,i,-N+y-1})^2z^2\D {\mathbb Q}(\xi_{i,y}<z)\\
\leq &\vare_n+\sum_{i=1}^n\sum_{y=1}^{2N+1}\int_{|z|>{\vare}|b_{n}|^{-1}}(b'_{n,i,-N+y-1})^2z^2\D {\mathbb Q}(\xi_{i,y}<z)\\
\leq&\vare_n+\int_{|z|>{\vare}|b_{n}|^{-1}}z^2\D {\mathbb Q}(\xi_{i,y}<z)\longrightarrow 0.
\end{align*}
It means the Lindeberg's condition satisfied and therefore $\frac{X_n}{A_n}\stackrel{D}{\longrightarrow}N(0,1)$. Since $A^2_n\sim \si^2n^H/\be^2$, we complete the proof of this proposition.\qed

\begin{remark}
{\rm This proposition is a CLT for weighted sum of stationary process. We use the method in \cite{Ibragimov1971Independent} to verify the Lindeberg's conditions holding. }\end{remark}

\begin{remark}
{\rm In the case of IID random environment, one has
\begin{align*}
&\EE_{\QQ}\bigg(\sum_{i=1}^n\sum_{x\in\ZZ}\om(i,x)\PP(S_i=x)\bigg)^2
=\sum_{i=1}^n2^{-2i}\binom{2i}{i}\\
=&\frac{2}{\pi}\sum_{i=1}^n\int_0^{\frac{\pi}{2}}(\sin x)^{2i}\D x
=\frac{2(2n+2)}{\pi}\int_0^{\frac{\pi}{2}}(\sin x)^{2n+2}\D x=(2n+2)2^{-(2n+2)}\binom{2n+2}{n+1}.
\end{align*}
An application of Stirling's formula shows that 
$$\EE_{\QQ}\bigg(\sum_{i=1}^n\sum_{x\in\ZZ}\om(i,x)\PP(S_i=x)\bigg)^2\sim \frac{2}{\sqrt{\pi}}\sqrt{n+1},
$$
whose order is smaller than the one in spatial correlated case. }
\end{remark}

\section{Gaussian fields, multiple stochastic integration and stochastic heat equations}
 
 In this section we will give a brief introduction about  Gaussian fields, multiple stochastic integrals, stochastic heat equation driven by time-white spatial-colored noise and its chaos expansion solution, for more details see \cite{Duncan2000Stochastic,Hu2001Heat,Hu2002Chaos}.
 
 \subsection{Gaussian fields, multiple stochastic integral.} 
 
 A spatially homogeneous Gaussian field that is white in time and correlated in space
is a mean zero Gaussian process $\{W(\phi), \phi\in\mathscr{S}([0,1]\times\RR)\}$,  defined on some probability space $(\Om',{\mathcal F}',{\mathbb P}')$, with covariance 
\begin{align}\label{colored field}
\EE'(W(\phi)W(\psi))=\int_{\RR^+}\int_{\RR}\Gamma(\D x)\phi(s,\cdot)\ast\tilde{\psi}(s,\cdot)(x)\D s,
\end{align}
where $\Gamma$ is a non-negative and non-negative definite tempered measure on $\RR$, $\ast$ is the convolution of two functions, $\tilde{\psi}(s,x)=\psi(s,-x)$ and $\mathscr{S}([0,1]\times\RR)$ is the space of rapidly decreasing functions on $[0,1]\times\RR$. 

In order to relate to the linear model introduced in the previous section, we are restricted ourselves in the paper the case that $\Ga$ is Riesz potential. That is, $\Gamma(\D x)=\frac{C_H\D x}{|x|^{2-2H}}=K(x)\D x$, with some $\frac{1}{2}<H<1$ and $C_H=H(2H-1)$, we may call $H$ Hurst parameter. $\Ga(\D x)$ has spectral measure $\mu(\D \xi)=\frac{H(2H-1)|\xi|^{1-2H}\D \xi}{2\Ga(2-2H)\cos(1-H)\pi}$ (see \eref{limitofspectrum}). It is obvious that when $H=\frac{1}{2}$ one has $\mu(\D \xi)=\frac{\D \xi}{2\pi}$, which corresponds to the white noise whose spectral density is constant. We denote the probability and expectation by $\PP_H$ and $\EE_H$, respectively. In this case, \eref{colored field} is
\begin{align}\label{fractional field}
\EE_{H}(W(\phi)W(\psi))=C_H\int_{0}^1\int_{\RR^2}\phi(s,x)|x-y|^{2H-2}\psi(s,y)\D x\D y\D s.
\end{align}

We introduce the following Hilbert space:
$$
\mathcal{L}_H=\{f:[0,1]\times{\mathbb R}\longrightarrow\RR; \|f\|_H^2=\int_0^1\int_{\mathbb R}\int_{\mathbb R} f(s,u)K(u,v)f(s,v)\D s\D u\D v<\infty\}, 
$$
where $K(u,v)=K(u-v)$. By Hardy-Littlewood inequalities \cite[Theorem 2.1]{Memina2001In} or \cite[Theorem 4.3]{Lieb1996}
we have, for some positive constant $A_H$,
\begin{align}\label{hl}
\left|\int_0^1\int_{\mathbb R}\int_{\mathbb R} f(s,u)K(u,v)f(s,v)\D s\D u\D v\right|\leq A_H\int_0^1\left(\int_{\mathbb R}|f(s,u)|^{\frac{1}{H}}\D u\right)^{2H}\D s.
\end{align}

For $t\in[0,1], A\in \mathcal{B}(\RR)$, choose an approximation sequence $f_n\in\mathcal{L}_H$ such that $f_n\downarrow 1_{[0,t]\times A}$, and define a worthy martingale measure in the sense of Walsh by
$$
W_t(A)\stackrel{\De}{=}W((0,t])\times A)=\lim_{n\rightarrow\infty}W(f_n).
$$
Hence
$$
\EE_H(W_t(A)W_s(B))=C_H(s\wedge t)\int_{A\times B}|x-y|^{2H-2}\D x\D y.
$$
Especially, we have a mean zero Gaussian random field $W=\{W(t,x)=W_t([0,x]):(t,x)\in[0,1]\times \RR\}$ with covariance
$$
\EE_H(W(t,x)W(s,y))=\frac{1}{2}(t\wedge s)||x|^{2H}+|y|^{2H}-|x-y|^{2H}|.
$$
We call $W$ fractional Gaussian field although the phrase ``fractional" is implicitly meant to time.
Taking the formal derivative $\dot{W}(t,x)$ of the field $W(t,x)$, i.e., $\dot{W}(t,x)=\frac{\partial^2W(t,x)}{\partial t\partial x}$, one can write the Gaussian process $W(f)$ as the following stochastic integral
\begin{equation}\label{noise}
W(f)=\int_0^1\int_{\RR}f(t,x)W(\D t\D x)=\int_0^1\int_{\mathbb R} f(t,x)\dot{W}(t,x)\D t\D x.
\end{equation}
for $f\in\mathcal{L}_H$.

Let $H_n$ be the Hermite polynomial of degree $n$, i.e.,
$$
H_n(x)=(-1)^n\e^{x^2/2}\frac{\D ^n}{\D x^n}\e^{-x^2/2},\quad x\in{\mathbb R}.
$$

We will use these Hermite polynomials to define multiple It\^{o} integrals \cite{Huang1997Introduction}. Let $\{h_1,h_2,\dots\}$ be an orthonormal basis of ${\mathcal L}_H$. Then $\{h_{i_1}\otimes h_{i_2}\otimes\dots\otimes h_{i_k}, i_1,i_2,\dots,i_k, k\in\NN\}$ is an orthonormal basis of ${\mathcal L}_H^{\otimes k}$, by which we denote the symmetric tensor product of ${\mathcal L}_H$. Here we still use the notation $\otimes$ to represent symmetric product instead of $\hat{\otimes}$.  Actually, we have
\begin{align}\label{ksym}
&{\mathcal L}_H^{\otimes k}=\{f: ([0,1]\times {\mathbb R})^k\rightarrow {\mathbb R}~\text{symmetric};\nonumber \\
&\int_{[0,1]^k}\int_{{\mathbb R}^{2k}}f(t_1,x_1,t_2,x_2,\dots,t_k,x_k)\prod_{i=1}^kK(x_i,y_i)f(t_1,y_1,t_2,y_2,\dots,t_k,y_k)\D {\mathbf t}\D {\mathbf x}\D {\mathbf y}<\infty\}
\end{align}
with ${\mathbf t}=(t_1,t_2,\dots,t_k),{\mathbf x}=(x_1,x_2,\dots,x_k)$ and ${\mathbf y}=(y_1,y_2,\dots,y_k)$.
For $f\in{\mathcal L}_H$, $f^{\otimes k}\in {\mathcal L}_H^{\otimes k}$, furthermore, if $\|f\|_H=1$, we define the multiple Ito integral of $f^{\otimes k}$ with respect to $\dot{W}$ by 
$$
\int_{([0,1]\times {\mathbb R})^k}f^{\otimes k}({\mathbf t},{\mathbf x})W^{\otimes k}(\D {\mathbf t}\D {\mathbf x})=H_k(W(f)):=I_k(f^{\otimes k}).
$$
Then the polarization procedure can be used to define the multiple integral of the form:
$$
\int_{([0,1]\times {\mathbb R})^k}f_1\otimes f_2\otimes\dots\otimes f_k({\mathbf t},{\mathbf x})W^{\otimes k}(\D {\mathbf t}\D {\mathbf x}):=I_k(f_1\otimes f_2\otimes\dots\otimes f_k).
$$
And we can go as usual to define $k$-multiple stochastic integral $I_k(f)$, 

$$
I_k(f)=\int_{([0,1]\times {\mathbb R})^k}f({\mathbf t},{\mathbf x})W^{\otimes k}(\D {\mathbf t}\D {\mathbf x}),
$$
for general symmetric functions $f$ in ${\mathcal L}_H^{\otimes k}$ by the density argument. Furthermore we have
$$
{\mathbb E}(I_k(f)I_k(g))=k!\langle f,g\rangle_{{\mathcal L}_H^{\otimes k}}.
$$
For $f\in {\mathcal L}_H^{\otimes m}, g\in{\mathcal L}_H^{\otimes n}, 1\leq r\leq m\land n$,  we define $r$-order contraction of two symmetry functions $f$ and $g$ by
\begin{align*}
&f\otimes_rg(t_1,x_1;\dots;t_{m+n-2r},x_{m+n-2r})\nonumber\\
=&{\text Sym}\left\{\int_{[0,1]^r}\int_{{\mathbb R}^{2r}}f(t_1,x_1;\dots;t_{n-r},x_{n-r};s_1,u_1;\dots;s_r,u_r)\right.\nonumber\\
&\times \Pi_{i=1}^rK(u_i,v_i)g(t_1,x_1;\dots;t_{m-r},x_{m-r};s_1,v_1;\dots;s_r,v_r)\D {\mathbf s}\D {\mathbf u}\D {\mathbf v}\bigg\},
\end{align*}
where ${\text Sym}\{\cdot\}$ means symmetrizing the arguments. Then one has the following recursive identities:
\begin{align}\label{contract1}
I_n(f)I_m(g)=\sum_{r=0}^{m\land n}r!\binom{n}{r}\binom{m}{r}I_{n+m-2r}(f\otimes_rg)
\end{align}
for $f\in {\mathcal L}_H^{\otimes m}, g\in{\mathcal L}_H^{\otimes n}$. Especially, when $m=1$, it is reduced to 
\begin{align}\label{contract2}
I_n(f)I_1(g)=I_{n+1}(f\otimes g)+nI_{n-1}(f\otimes_1g).
\end{align}

For later use, we give an example of the contraction of two functions. Assume $f_1,f_2\in {\mathcal L}_H, m,n\in{\mathbb N}$. Notice that the symmetrical function $f_1^{\otimes m}\otimes f_2^{\otimes(n-1)}$ has $\binom{m+n-1}{m}$ terms, where $\binom{m+n-2}{m-1}$ terms end with $f_1$ and $\binom{m+n-2}{n-2}$ terms end with $f_2$. Hence we have 
\begin{align*}
&(f_1^{\otimes m}\otimes f_2^{\otimes(n-1)})\otimes_1f_2\\
=&\frac{m}{m+n-1}f_1^{\otimes (m-1)}\otimes f_2^{\otimes(n-1)}\langle f_1,f_2\rangle_{H}\\
    &+\frac{n-1}{m+n-1}f_1^{\otimes m}\otimes f_2^{\otimes(n-2)}\|f_2\|^2_{H}.
\end{align*}
Furthermore, according to \eqref{contract2}, we have
\begin{align}\label{contract3}
&I_{m+n}(f_1^{\otimes m}\otimes f_2^{\otimes n})\nonumber\\
=&I_{m+n-1}(f_1^{\otimes m}\otimes f_2^{\otimes(n-1)})I_1(f_2)-m\langle f_1,f_2\rangle _HI_{m+n-2}(f_1^{\otimes (m-1)}\otimes f_2^{\otimes(n-1)})\\
&-(n-1)\|f_2\|^2_HI_{m+n-2}(f_1^{\otimes m}\otimes f_2^{\otimes(n-2)}).\nonumber
\end{align}
      
      Now we have the following chaos expansion results for square integrable variables.
\begin{prop}
Let $W$ be the gaussian random field above with spatial parameter $1/2<H<1$. Let $(\Om_H, {\mathcal F}_H, P_H)$ be the canonical probability space corresponding to $W$.  Then for any $F\in L^2(\Om_H)$, it admits the following chaos expansion:
$$
F=\sum_{k=0}^\infty I_k(f_k),
$$
where $f_k\in {\mathcal L}_H^{\otimes k}, k=0,1,\dots,$ and the series converges in $L^2(\Om_H,{\mathcal F}_H,P_H)$. Moreover,
$$
{\mathbb E}_H[F^2]=\sum_{k=0}^\infty {k!}\|f_k\|^2_{{\mathcal L}_H^{\otimes k}}.
$$
\end{prop}

\begin{remark}\label{remmultipleintegral}
{\rm (i) It is obvious that $I_k(h^{\otimes k})\ne [I_1(h)]^k$, and $I_k(h_1\otimes h_2\otimes\dots\otimes h_k)\ne\prod_{i=1}^kI_1(h_i)$ unless $h_1,h_2,\dots,h_k$ are orthogonal in ${\mathcal L}_H$. Usually, $I_k(h^{\otimes k})$ is called $k$-order Wick product of $I_1(h)$. 

(ii) Denote by $H(W)$ the linear space spanned by $\{W(t,x); (t,x)\in {\mathbb R}^+\times{\mathbb R}\}$, then there is a unique isomorphism $\Phi$ between $\oplus_{k=0}^\infty [H(W)]^{\otimes k}$ and $L^2(\Om_H,{\mathcal F}_H,P_H)$ such that if $\{\zeta_i:i\in\NN\}$ is a complete orthonormal system (CONS) in $H(W)$, then the family 
$$
\left(\frac{k!}{k_{\ga_1}!\cdots k_{\ga_i}!}\right)^{\frac{1}{2}}\Phi(\zeta_{\ga_1}^{\otimes k_{\ga_1}}\otimes\cdots\otimes\zeta_{\ga_i}^{\otimes k_{\ga_i}})=\prod_{j=1}^i(k_{\ga_j}!)^{-\frac{1}{2}}H_{k_{\ga_j}}(\zeta_{\ga_j}),
$$
$k\geq 0, i\geq 1, k_{\ga_1}+\cdots+k_{\ga_{i}}=k, \ga_1<\cdots<\ga_i$, is a CONS in $L^2(\Om_H,{\mathcal F}_H,P_H)$. See Huang, et al. \cite[P.595]{Huang1978Stochastic}. Hence, it suffices to define the multiple integral of the form $f_{\ga_1}^{\otimes k_{\ga_1}}\otimes\cdots\otimes f_{\ga_i}^{\otimes k_{\ga_i}}$ in ${\mathcal L}_H^{\otimes k}$.}
\end{remark}

\subsection {Heat equations with white noise potentials.} We turn to stochastic heat equations \eqref{SHE1} with multiplicative noise 
and initial value $u(s,x)=u(x)$, $0\leq s\leq t\leq 1, x\in{\mathbb R}$. Its solution is formulated in the mild form, i.e.,
\begin{equation}\label{mildsolution}
u(t,x;s)=P_{t-s}u(x)+\sqrt{2}\be\int_s^t\int_{\mathbb R}P({t-r},x-z)u(r,z) W(\D r\D z),
\end{equation}
where $P(t,x)=\frac{1}{\sqrt{2\pi t}}\e^{-\frac{x^2}{2t}}$ and $P_tf(x)=\int_{\mathbb R}\frac{1}{\sqrt{2\pi t}}\e^{-\frac{(x-y)^2}{2t}}f(y)\D y$. Furthermore, let $u(x)=\de(x-y)$, the Dirac delta function at zero, we get a four-parameter field $u(t,x;s,y)$ by
$$
u(t,x;s,y)=P({t-s},x-y)+\sqrt{2}\be\int_s^t\int_{\mathbb R}P({t-r},x-z)u(r,z;s,y)W(\D r\D z).
$$
Iterating the equation yields a formal chaos expansion for $u(t,x;s,y)$:
\begin{align}\label{chaossolution}
&u(t,x;s,y)\nonumber\\
=&P({t-s},x-y)\nonumber\\
&+\sum_{k=1}^\infty(\sqrt{2}\be)^k\int_{\De(s,t]^k}\int_{{\mathbb R}^k}
\Pi_{i=1}^kP(t_i-t_{i-1},x_i-x_{i-1})P({t-t_{k}},x-x_{k})W(\D t_i\D x_i)\\
=&P_{t-s}(x-y)
+\sum_{k=1}^\infty(\sqrt{2}\be)^kI_k(\widetilde{P_k(t,x;s,y)})\nonumber
\end{align} 
with $\De(s,t]^k=\{s<t_1<\cdots<t_k<t\}$. Put
\begin{align}\label{k-td}
P_k(t,x;s,y;t_1,\dots,t_k;x_1,\dots,x_k)=&
\Pi_{i=1}^{k+1}P_{t_i-t_{i-1}}(x_i-x_{i-1})\stackrel{\De}{=}P_k(s,y;t,x;{\mathbf t},{\mathbf x})
\end{align}
(or $P_k(t,x;{\mathbf t},{\mathbf x})$ when $s=0, y=0$) is the transition density function of Brownian motion from $y$ to $x$ through $k$ points on the time interval $(s,t]$, where ${\mathbf t}=(t_1,\cdots,t_k), {\mathbf x}=(x_1,\cdots,x_k)$,  $t_0=s, x_0=y, t_{k+1}=t,x_{k+1}=x$, and $\widetilde{P_k(t,x;s,y)}$ is the symmetrization of $P_k(t,x;s,y;\cdot)$ in the variables $(t_1,s_1)$, $\cdots, (t_k,s_k)$.

\begin{prop}
For any $s\geq 0, y\in\RR$, $u(t,x;s,y)$ given by \eqref{chaossolution} is the unique solution to equation \eqref{SHE1} with initial data $u(s,x;s,y)=\de(x-y)$.
\end{prop}
\proof We compute the ${L}^2$ norm of each chaos. Denote
$$
\Theta_k(t,x;s,y)=k!(\sqrt{2}\be)^{2k}{\mathbb E}_H(I_k(\widetilde{P_k(t,x;s,y)}))^2,\quad\text{for}~k=1,2,\dots.
$$
Then by the isometric equality, we have
\begin{align*}
\Theta_k(t,x;s,y)=(\sqrt{2}\be)^{2k}\int_{\De(s,t]^k}\int_{{\mathbb R}^{2k}}\Pi_{i=1}^kK(x_i,y_i)P_k(t,x;s,y;\tau;{\mathbf x})P_k(t,x;s,y;\tau;{\mathbf y})\D {\mathbf x}\D {\mathbf y}\D \tau,
\end{align*}
where ${\mathbf x}=(x_1,\dots,x_k),{\mathbf y}=(y_1,\dots,y_k)$. 

By the proof of Lemma 6.1 in \cite{Hu2001Heat}, we know
\begin{align*}
\int_{{\mathbb R}^2}P({t_i-t_{i-1}},x_i-x_{i-1})K(x_i,y_i)P({t_i-t_{i-1}},y_i-y_{i-1})\D x_{i}\D y_{i}\leq A_H(t_i-t_{i-1})^{H-1}
\end{align*}
for $i=1,2,\dots,k$, where $A_H$ is a generic constant depending only on $H$. Hence, by Cauchy-Schwarz inequality, we have
\begin{align*}
&\int_{{\mathbb R}^2}P({t_1-s},x_1-y)P({t_2-t_1},x_2-x_1)K(x_1,y_1)P({t_1-s},y_1-y)P({t_2-t_1},y_2-y_1)\D x_{1}\D y_{1}\\
\leq&\frac{1}{2\pi\sqrt{(t_1-s)(t_2-t_1)}}\left(\int_{\mathbb R}P({t_1-s},\sqrt{2}(x_1-y))K(x_1,y_1)P({t_1-s},\sqrt{2}(y_1-y))\D x_{1}\D y_{1}\right)^{\frac{1}{2}}\\
&~~~~~~~~~~~~~~~~~~~~~~~~~~\times \left(\int_{\mathbb R}P({t_2-t_1},\sqrt{2}(x_2-x_1))K(x_1,y_1)P({t_2-t_1},\sqrt{2}(y_2-y_1))\D x_{1}\D y_{1}\right)^{\frac{1}{2}}\\
\leq&\frac{A_H}{2^{1+H}\pi}(t_1-s)^{\frac{H}{2}-1}(t_2-t_1)^{\frac{H}{2}-1}.
\end{align*}

It follows that
\begin{align*}
\Theta_k(t,x;s,y)&\leq\frac{A_H^k(\sqrt{2}\be)^{2k}}{2^{1+H}\pi}\int_{\De(s,t]^k}(t_1-s)^{\frac{H}{2}-1}(t_2-t_1)^{\frac{H}{2}-1}\Pi_{i=3}^k(t_i-t_{i-1})^{H-1}\D t_1\dots\D t_k\\&=\frac{A^k_H2^k\be^{2k}\Ga^2(\frac{H}{2})\Ga^{k-1}(H)}{2^{(1+H)}\pi \Ga(kH)}(t-s)^{kH-1},
\end{align*}
and 
$$
{\mathbb E}_H(u(t,x;s,y))^2\leq P^2({t-s},x-y)+\sum_{k=1}^\infty\Theta_k(t,x;s,y)<\infty.
$$
Therefore, the chaos expansion \eqref{chaossolution} is the unique solution to equation \eqref{SHE1} in $ L^2$.\qed

\section{The convergence of partition functions in the sense of finite dimensional distribution} 

In this section we are focused on the convergence of partition function for the polymer measure in the sense of finite dimensional distribution. Firstly, we consider the convergence for modified partition function. 

\subsection {Modified partition function.} The scaled modified function $\mathfrak {Z}_n^\om(\be n^{-\vr})$ and its expansion are given by
\begin{align}
\label{mpartition12}
&\mathfrak {Z}_n^\om(\be n^{-\vr})={\mathbb E_{\mathbb P}}\left[\prod_{i=1}^n(1+\be n^{-\vr}\om(i,S_i))\right],
\end{align}
and
\begin{align}
\label{mpartition13}
\mathfrak {Z}_n^\om(\be n^{-\vr})
=&{\mathbb E_{\mathbb P}}\left[1+\sum_{k=1}^n\be^kn^{-k\vr}\sum_{{\mathbf i}\in D_k^n}\prod_{j=1}^k\om(i_j,S_{i_j})\right]\nonumber\\
=&1+\sum_{k=1}^n\be^kn^{-k\vr}\sum_{{\mathbf i}\in D_k^n}\sum_{{\mathbf x}\in {\mathbb Z}^k}\left[\prod_{j=1}^k\om(i_j,x_j)p({i_j-i_{j-1}},{x_j-x_{j-1}})\right]\\
=&1+\sum_{k=1}^n\be^kn^{-k\vr}\sum_{{\mathbf i}\in D_k^n}\sum_{{\mathbf x}\in {\mathbb Z}^k}\om({\mathbf i},{\mathbf x})p_k({\mathbf i},{\mathbf x})\nonumber.
\end{align}
See \eqref{porder1} for $D_k^n$ and \eqref{k-transition} for $p_k({\mathbf i},{\mathbf x})$. This is a point-to-line (starting from zero) modified partition function. We have also point-to-point modified partition function as follows. For ${\mathbb Z}\ni x\leftrightarrow n, 0\leq m<n, y\in{\mathbb Z}$
\begin{align*}
&\mathfrak {Z}_n^\om(m,y;n,x;\be)\\
=&{\mathbb E_{\mathbb P}}\left[\prod_{i=m+1}^n(1+\be\om(i,S_i)){\mathbbm 1}_{\{S_n=x\}}|_{S_m=y}\right]\\
=&1+{\mathbb E_{\mathbb P}}\left[\sum_{k=1}^{n-m}\be^k\sum_{{\mathbf i}-m\in D_k^{n-m}}\left.\prod_{j=1}^k\om(i_j,S_{i_j})\right|S_m=y, S_n=x\right]p(n-m,x-y)\\
=&1+\sum_{k=1}^{n-m}\be^k\sum_{{\mathbf i}-m\in D_k^{n-m}}\sum_{{\mathbf x}\in {\mathbb Z}^k}\left[\prod_{j=1}^k\om(i_j,x_j)p_{y,x}^{n-m}({\mathbf i},{\mathbf x})\right]p(n-m,x-y),
\end{align*}
where $p_{y,x}^{n-m}({\mathbf i},{\mathbf x})=\frac{p(n-m-i_k,x-x_k)}{p(n-m,x-y)}\prod_{j=1}^kp(i_j-i_{j-1},x_j-x_{j-1})$ is the transition kernel for random walks conditioned to be at position $x$ at time $n$, and $p(n-m,x-y)=\PP(S_{n-m}=x|S_0=y)$. When $m=0$ and $y=0$, we denote it by $\mathfrak {Z}_n^\om(n,x;\be)$. Similarly, we have ${Z}_n^\om(m,y;n,x;\be)$ and ${Z}_n^\om(n,x;\be)$ as \eqref{partition2}.

We are aimed at the asymptotic behavior of \eqref{mpartition13}, but begin with general $f$ instead of the transition density $p$. For $k\leq n\in {\mathbb N}$ we use the notations ${\mathfrak R}_k^n$ as in \cite{alberts2014} to denote the set of rectangles of the form 
$$
{\mathfrak R}_k^n\stackrel{\De}{=}\left\{\left(\frac{{\mathbf i}-1}{n},\frac{{\mathbf i}}{n}\right]\times\left(\frac{{\mathbf x}-1}{\sqrt{n}},\frac{{\mathbf x}+1}{\sqrt{n}}\right]:{\mathbf i}\in E_k^n, {\mathbf i}\leftrightarrow{\mathbf x}\right\},
$$
where
\begin{align}\label {norder1}
E_k^n=\{{\textbf i}=(i_1,i_2,\dots,i_k)\in [n]^k: 1\leq i_j\ne i_l\leq n, {\text {for~}} j\ne l\leq k\}.
\end{align}

For $f\in L^2([0,1]^k\times {\mathbb R}^k)$ define $\bar{f}_n$ by
$$
{\bar f}_n({\mathbf t},{\mathbf x})=\frac{1}{|R|}\int_{R}f\D {\mathbf t}\D{\mathbf x},\quad ({\mathbf t},{\mathbf x})\in R\in{\mathfrak R}_k^n,
$$
and weighted $U-$statistics ${\mathcal S}_k^n$ by 
\begin{align}\label{Uk}
{\mathcal S}_k^n(f)=2^{k/2}\sum_{{\mathbf i}\in E_k^n}\sum_{{\mathbf x}\in {\mathbb Z}^k}
{\bar f}_n\big(\frac{{\mathbf i}}{n},\frac{{\mathbf x}}{\sqrt{n}}\big)\om({\mathbf i},{\mathbf x})\mathbbm{1}_{\{{\mathbf i}\leftrightarrow{\mathbf x}\}},
\end{align}
$|R|$ is the Lebesgue measure of $R$.  Here the introduction of factor $\sqrt{2}$ before the sum is to cancel the parity of random walk.

Notice that $p_k^n$ is a constant on each rectangle in ${\mathfrak R}_k^n$, so that $\bar{p}_k^n=p_k^n$. Also we have, for ${\mathbf i}\in D_k^n,{\mathbf x}\in{\mathbb Z}^k$ with ${\mathbf i}\leftrightarrow{\mathbf x}$,
$$
p_k^n\left(\frac{{\mathbf i}}{n},\frac{{\mathbf x}}{\sqrt{n}}\right)=\bar{p}_k^n=2^{-k}p_k({\mathbf i}, {\mathbf x}).
$$
Then, we have 
\begin{align}\label{U2}
{\mathcal S}_k^n(p_k^n)=2^{k/2}\sum_{{\mathbf i}\in E_k^n}\sum_{{\mathbf x}\in {\mathbb Z}^k}
{\bar f}_n\big(\frac{{\mathbf i}}{n},\frac{{\mathbf x}}{\sqrt{n}}\big)\om({\mathbf i},{\mathbf x})\mathbbm{1}_{\{{\mathbf i}\leftrightarrow{\mathbf x}\}}=2^{-k/2}\sum_{{\mathbf i}\in D_k^n}\sum_{{\mathbf x}\in {\mathbb Z}^k}
\om({\mathbf i},{\mathbf x})p_k({\mathbf i},{\mathbf x}),
\end{align}
which is the $k$th term in \eqref{mpartition13} except the coefficient. 

For general $f\in L^2([0,1]^k\times \RR^k)$ with
$$
\int_{[0,1]^k}\int_{{\mathbb R}^{2k}}f(t_1,x_1,t_2,x_2,\dots,t_k,x_k)\prod_{i=1}^kK(x_i,y_i)f(t_1,y_1,t_2,y_2,\dots,t_k,y_k)\D {\mathbf t}\D {\mathbf x}\D {\mathbf y}<\infty,
$$
we have the following key point for obtaining the convergence of modified partition function \eqref{mpartition13}.

\begin{thm}\label{mainthm}
Let $f\in {\mathcal L}_H^{\otimes k}$. Then, as $n\rightarrow\infty$, 
\begin{align}\label{convergence1}
n^{-\frac{(H+1)k}{2}}{\mathcal S}_k^n(f)\stackrel{D}{\longrightarrow}\int_{[0,1]^k}\int_{\RR^k}
f({\mathbf t},{\mathbf x})W^{\otimes k}(\D{\mathbf t}\D{\mathbf x})=I_k(f),
\end{align}
where $W(\D t\D x)$ is fractional Gaussian noise \eqref{noise} and $I_k,k=1,2\dots,$ is $k$-multiple integral defined in Section 3. Furthermore, we have the following joint convergence
\begin{equation}\label{jointconvergence}
\bigg(n^{-\frac{(H+1)k_1}{2}}{\mathcal S}_{k_1}^n(f_1), \cdots,n^{-\frac{(H+1)k_r}{2}}{\mathcal S}_{k_r}^n(f_r)\bigg)\stackrel{D}{\longrightarrow}(I_{k_1}(f_1),\cdots,I_{k_r}(f_r))
\end{equation}
as $n\to\infty$ for $f_1\in{\mathcal L}_H^{\otimes k_1},\dots,f_r\in{\mathcal L}_H^{\otimes k_r}, r\in\NN$.
\end{thm}

The proof for theorem \ref{mainthm} is based on the following lemmas, whose proof will be postponed to the next subsection.

\begin{lemma}\label{chaosapproximation}
For all fixed $k,n$, ${\mathcal S}_k^n(f)$ is linear in $f$ with probability one, and for $k_1\ne k_2$, ${\mathbb E}_{\mathbb Q}({\mathcal S}_{k_1}^n(f_1){\mathcal S}_{k_2}^n(f_2))=0$ for $f_i\in{\mathcal L}_H^{\otimes k_i}, i=1,2$. Furthermore, for $k_1=k_2=k$, we have
\begin{align*}
{\mathbb E}_{{\mathbb Q}} [({\mathcal S}_{k}^n(f))^2]\leq C\la^kn^{(1+H)k}\|f\|^2_{H^k}
\end{align*}
for some generic positive constant $C$.
\end{lemma} 

\begin{lemma}\label{k=1}
The conclusion in Theorem \ref{mainthm} holds when $k=1$ for $f\in{\mathcal L}_{ H}$.
\end{lemma}

\begin{lemma}\label{k>1}
The conclusion in Theorem \ref{mainthm} holds for $f\in {\mathcal L}^k_H$ of the form $f=g^{\otimes k}$ with 
\begin{align}\label{indicator}
 g(t,x)=\mathbbm{1}_{[t_0,t_1]\times[x_0, x_1]}(t,x)\in{\mathcal L}_H
 \end{align}
 for $0\leq t_0\leq t_1\leq 1, x_0\leq x_1$ and $k>1$.
\end{lemma}

To implement linear extension to general $f$, according to Remark \ref{remmultipleintegral} (ii), we need the following result.

\begin{lemma}\label{tensor}
The conclusion in Theorem \ref{mainthm} holds for $f\in {\mathcal L}^k_H$ of the form $f=g_1^{\otimes k_1}\otimes\cdots\otimes g_s^{\otimes k_s}$ with $g_1,\dots,g_s\in{\mathcal L}_H$ and $ k_1+\cdots+k_s=k, k_1>0,\dots, k_s>0, s=2,3,\cdots$. 
\end{lemma}

{\it Proof of Theorem \ref{mainthm}}\quad Notice that ${\mathcal S}_k^n(f)$ has symmetrizing property, we assume that $f$ is symmetrical. Let $\{h_i\}_{i=1}^\infty$ be a complete orthonormal base of ${\mathcal L}_H$. By lemma \ref{tensor}, we know the conclusion of theorem \ref{mainthm} is true for all $f$ of the form $f=h_{i_1}^{\otimes k_{i_1}}\otimes\cdots\otimes h_{i_s}^{\otimes k_{i_s}}$ with $k_{i_1}>0,\cdots, k_{i_s}>0$ and $k_{i_1}+\cdots+k_{i_s}=k$. Furthermore, it is also true for the linear combination of such f's, since by the recursive identity \eqref{contract2} and the proofs of Lemma \ref{k>1} and Lemma \ref{tensor} all things as well as the joint convergence are reduced to the case of the convergence of ${\mathcal S}_1^n(f)$ by continuous mapping theorem. \qed

The next two results state the convergence of modified partition function. 

\begin{thm}\label{modcon}
Assume $\om$ given by \eqref{envir1}. Then there exists a square integrable random variable $\rZ_{\sqrt{2}\be}$, such that 
$$
\mathfrak {Z}^\om_n(\be n^{-\vr})\stackrel{D}{\longrightarrow}\rZ_{\sqrt{2}\be},
$$
as $n\to\infty$. Actually, $\rZ_{\sqrt{2}\be}$ has a chaos decomposition. 
\end{thm}
\proof Due to \eqref{U2}, we have 
\begin{equation}\label{mpconvergence}
\begin{split}
\mathfrak {Z}_n^\om(\be n^{-\vr})=&1+\sum_{k=1}^n\be^kn^{-k\vr}\sum_{{\mathbf i}\in D_k^n}\sum_{{\mathbf x}\in {\mathbb Z}^k}\bigg[\prod_{j=1}^k\om(i_j,x_j)p_k({\mathbf i},{\mathbf x})\bigg]\\
=&1+\sum_{k=1}^n2^{k/2}\be^kn^{-\frac{(H+1)k}{2}}{\mathcal S}_k^n(n^{\frac{k}{2}}p^n_k).
\end{split}
\end{equation}
As the proof of Lemma 4.4 in \cite{alberts2014}, one has 
\begin{equation}\label{C1}
I_n(P_{\sqrt{2}\be}):=1+\sum_{k=1}^\infty n^{-\frac{(H+1)k}{2}}{\mathcal S}_k^n((\sqrt{2}\be)^kP_k )\longrightarrow I(P_{\sqrt{2}\be}):=1+\sum_{k=1}^\infty I_k((\sqrt{2}\be)^kP_k)
\end{equation}
as $n\to\infty$, where $P_{\sqrt{2}\be}=(1,\sqrt{2}\be P_1,\cdots, (\sqrt{2}\be)^kP_k,\cdots)\in\oplus_{k=0}^\infty{\mathcal L}_H^{\otimes k}$ and for $k\geq 1$, $P_k(\tau,{\bf x})=\Pi_{i=1}^kp(\tau_i-\tau_{i-1},x_i-x_{i-1}), \tau=(\tau_1,\cdots,\tau_k), {\bf x}=(x_1,\cdots,x_k)$ with $\tau_0=0, x_0=0$; $P_0=1.$ This is based on the following basic facts. By the isometry property of  $I$, one has $\operatorname { that } \sum _ { k = 0 } ^ { M } I _ { k } \left( g _ { k } \right) \rightarrow \sum _ { k = 0 } ^ { \infty } I _ { k } \left( g _ { k } \right) \operatorname { in } L ^ { 2 } \left( \Omega _ { H } , \mathcal { F } _ { H } , \mathbb { P }_H \right)$ as $M \rightarrow \infty$. Also one has
$$
\sum _ { k = 0 } ^ { M } n^{-\frac{(H+1)k}{2}} \mathcal { S } _ { k } ^ { n } \left((\sqrt{2}\be)^kP_k \right) \stackrel { M \rightarrow \infty } { \longrightarrow } \sum _ { k = 0 } ^ { \infty } n^{-\frac{(H+1)k}{2}}\mathcal { S } _ { k } ^ { n } \left((\sqrt{2}\be)^kP_k \right)
$$
in $L^2(\QQ)$, uniformly in $n$, by Lemma \ref{chaosapproximation}. Combining \eqref{convergence1} and \eqref{jointconvergence} yields
$$
\sum _ { k = 0 } ^ { M } n^{-\frac{(H+1)k}{2}} \mathcal { S } _ { k } ^ { n } \left( (\sqrt{2}\be)^kP_k  \right) \stackrel { D } { \longrightarrow } \sum _ { k = 0 } ^ { M } I _ { k } \left( (\sqrt{2}\be)^kP_k  \right).
$$
Lemma 4.2 in \cite{alberts2014} (also see \cite[Chapter 1, Theorem 4.2]{Billi}) shows \eqref{C1} true.
Now by considering the difference between $\mathfrak {Z}_n^\om(\be n^{-\vr})$ and $I_n(P_{\sqrt{2}\be})$ one has
$$
\begin{aligned}
&\mathfrak {Z}_n^\om(\be n^{-\vr})-I_n(P_{\sqrt{2}\be})\\
=&\sum _ { k = 0 } ^ { n } (\sqrt{2} \beta) ^ { k } n^{-\frac{(H+1)k}{2}}\mathcal { S } _ { k } ^ { n } ( P _ { k } - n ^ { k / 2 } p _ { k } ^ { n } ) + \sum _ { k = n + 1 }^\infty (\sqrt{2} \beta) ^ { k } n^{-\frac{(H+1)k}{2}} \mathcal { S } _ { k } ^ { n } ( P_ { k } ).
\end{aligned}
$$
 By using \eqref{hl}, which says the norm of $\mathcal {L}_H$ is controlled by $\mathcal{L}_2$ norm, combining local limit theorem and the proof of Proposition 5.3 in \cite{alberts2014} gives the above quantity converging to zero, thus 
$$
\mathfrak {Z}_n^\om(\be n^{-\vr})\longrightarrow I(P_{\sqrt{2}\be})=1+\sum_{k=1}^\infty I_k((\sqrt{2}\be)^k P_k):=\rZ_{\sqrt{2}\be}
$$
 as $n\to\infty$.
\qed

\begin{thm}\label{convergence}
Let $\{u(t,x),(t,x)\in [0,1]\times{\mathbb R}\}$ be the solution to \eqref{SHE1} with initial data $u(x)=\de(x)$. Then
\begin{align*}
\frac{\sqrt{n}}{2}\mathfrak {Z}^\om(nt,\sqrt{n}x;\be n^{-\vr})\longrightarrow u(t,x) \quad \mbox{weakly, as}\quad n\longrightarrow\infty.
\end{align*}
\end{thm}
\proof As before, we have 
\begin{align*}
&\mathfrak {Z}_n^\om(n,x;\be)\\
=&\EE_{\PP}\bigg(\Pi_{i=1}^n(1+\be\om(i,S_i))\mathbbm{1}_{\{S_n=x\}}\bigg)\\
=&\sum_{k=0}^n\be^k\sum_{{\mathbf i}\in D_k^n}\sum_{{\mathbf x}\in \ZZ^k}\om({\mathbf i},{\mathbf x})\Pi_{j=0}^kp(i_{j+1}-i_j,x_{j+1}-x_j),
\end{align*}
with $i_{k+1}=n,x_{k+1}=x,i_0=0$ and $x_{i_0}=0$. By recalling the definition of $p_k^n$ (see \eqref{SP1}), one has 
\begin{align*}
&\frac{\sqrt{n}}{2}\mathfrak {Z}_n^\om(nt,\sqrt{n}x;\be n^{-\vr})\\
=&\frac{\sqrt{n}}{2}p_0^n(t,x)+\sum_{k=1}^n2^{k/2}\be^kn^{-\frac{(H+1)k}{2}}{\mathcal S}_k^n(n^{\frac{k+1}{2}}p^n_k(t,x;\cdot,\cdot)),
\end{align*}
where $p^n_k(t,x;{\mathbf t},\mathbf{x})=\frac{1}{2}p^n_k({\mathbf t},\mathbf{x})p^n_0(\lfloor nt\rfloor-\lfloor nt_k\rfloor, \sqrt{n}(x-x_k))$.
Then local limit theorem and the same proof as that of Theorem \ref{modcon} yield the convergence 
$$
\frac{\sqrt{n}}{2}\mathfrak {Z}_n^\om(nt,\sqrt{n}x;\be n^{-\vr})\longrightarrow P(t,x)+\sum_{k=1}^\infty I_k((\sqrt{2}\be)^k P_k(t,x;\cdot,\cdot))
$$
as $n\to\infty$. See \eqref{k-td} for the notation $P_k(t,x;{\mathbf t},{\mathbf x})$. Since 
$ P(t,x)+\sum_{k=1}^\infty I_k((\sqrt{2}\be)^k P_k(t,x;\cdot,\cdot))
$ is the solution to \eqref{SHE1} with initial data $u(x)=\de(x)$, we complete the proof of this theorem.\qed

\begin{remark}
{\rm (1) Similarly, one has the convergence of four-parameter field to the solution to equation \eqref{SHE1} with initial data $u(x)=\de(x-y), y\in{\mathbb R}$, i.e., 
$$
\frac{\sqrt{n}}{2}\mathfrak {Z}^\om_n(ns,\sqrt{n}y;nt,\sqrt{n}x;\be n^{-\vr})\stackrel{D}{\longrightarrow} u(s,y;t,x)
$$
as $n\rightarrow\infty$, where $u(s,y;t,x)$ is the solution to equation \eqref{SHE1} with initial data $u(s,y;s,x)=\de(x-y)$ for $y\in \RR$.

(2) The proof of tightness is deferred to Section 5.}
\end{remark}
\vskip 0.05in

\subsection {The proof for Lemmas.}  In this subsection we prove Lemmas 4.2-5.

{\it Proof of Lemma \ref{chaosapproximation}}.\quad We only show the last claim. For $f\in{\mathcal L}_H^{\otimes k}$
\begin{align*}
&{\mathbb E}_{{\mathbb Q}} [({\mathcal S}_{k}^n(f))^2]\\
=&2^k\sum_{{\mathbf i}\in E_k^n}\sum_{{\mathbf x}\in {\mathcal E}_i^k}\sum_{{\mathbf y}\in {\mathcal E}_i^k}{\bar f}_n(\frac{{\mathbf i}}{n},\frac{{\mathbf x}}{\sqrt{n}}){\bar f}_n(\frac{{\mathbf i}}{n},\frac{{\mathbf y}}{\sqrt{n}})\ga(|{\mathbf x}-{\mathbf y}|)\mathbbm{1}_{\{{\mathbf i}\leftrightarrow{\mathbf x}\}}1_{\{{\mathbf i}\leftrightarrow{\mathbf y}\}}\\
\leq &\sum_{{\mathbf i}\in E_k^n}\sum_{{\mathbf x}\in {\mathcal E}_i^k}\sum_{{\mathbf y}\in {\mathcal E}_i^k}{\bar f}_n(\frac{{\mathbf i}}{n},\frac{{\mathbf x}}{\sqrt{n}}){\bar f}_n(\frac{{\mathbf i}}{n},\frac{{\mathbf y}}{\sqrt{n}})\ga(|{\mathbf x}-{\mathbf y}|)\\
\leq &C\la^kn^{(1+H)k}\sum_{{\mathbf i}\in E_k^n}\sum_{{\mathbf x}\in {\mathcal E}_i^k}\sum_{{\mathbf y}\in {\mathcal E}_i^k}{\bar f}_n(\frac{{\mathbf i}}{n},\frac{{\mathbf x}}{\sqrt{n}}){\bar f}_n(\frac{{\mathbf i}}{n},\frac{{\mathbf y}}{\sqrt{n}})\left|\frac{{\mathbf x}-{\mathbf y}}{\sqrt{n}}\right|^{1-2\al}\frac{1}{n^k}\frac{1}{{\sqrt{n}^{2k}}}\\
=&C\la^kn^{(1+H)k}\sum_{{\mathbf i}\in E_k^n}\sum_{{\mathbf x}\in {\mathcal E}_i^k}\sum_{{\mathbf y}\in {\mathcal E}_i^k}{\bar f}_n(\frac{{\mathbf i}}{n},\frac{{\mathbf x}}{\sqrt{n}}){\bar f}_n(\frac{{\mathbf i}}{n},\frac{{\mathbf y}}{\sqrt{n}})\left|\frac{{\mathbf x}-{\mathbf y}}{\sqrt{n}}\right|^{1-2\al}\frac{1}{n^k}\frac{1}{{\sqrt{n}^{2k}}}\\
\leq&C\la^kn^{(1+H)k}\int_{[0,1]^k}\int_{\RR^{2k}}|{\bar f}_n({\mathbf t},{\mathbf x})||{\bar f}_n({\mathbf t},{\mathbf y})|\left|{\mathbf x}-{\mathbf y}\right|^{1-2\al}\D{\mathbf t}\D{\mathbf x}\D{\mathbf y}\\
=&Cn^{(1+H)k}\|{\bar f}_n\|^2_{{H}^k}\leq Cn^{(1+H)k}\|f\|^2_{{H}^k}.
\end{align*}\qed

{\it Proof of Lemma \ref{k=1}}.\quad  First, we show it is true for $f=\mathbbm{1}_{[t_0,t_1]\times[x_0, x_1]}$ for some $0\leq t_0<t_1\leq 1, x_0<x_1.$
 In this case, 
 $$
 {\mathcal S}_1^n(f)=2^{1/2}\sum_{\{ nt_0\leq { i}\leq nt_1\}}\sum_{\{\sqrt{n}x_0\leq {x}\leq\sqrt{n}x_1\}}
\om({i},{ x}){\mathbbm1}_{\{{ i}\leftrightarrow{x}\}}
$$
and
\begin{align*}
&{\mathbb E}_{{\mathbb Q}} [({\mathcal S}_{1}^n(f))^2]=n^{\frac{3}{2}}(t_1-t_0)(x_1-x_0)\ga(0)+2n(t_1-t_0)\sum_{k=1}^{N-1}(N-k)\ga(k)+O(n),
\end{align*}
here and hereafter the sums for $i$ and $x$ are over $E_k^n$ and $\ZZ$, respectively, $N=\lfloor\sqrt{n}(x_1-x_0)\rfloor$ is the largest integer smaller than $\sqrt{n}(x_1-x_0)$, unless other statements.
Notice that the second term of ${\mathbb E}_{{\mathbb Q}} [({\mathcal S}_{1}^n(f))^2]$ can be written as
\begin{align*}
&2\la n(t_1-t_0)N^{3-2\al}\sum_{k=1}^{N-1}(1-\frac{k}{N})\left(\frac{k}{N}\right)^{1-2\al}\frac{1}{N}\\
\sim&2\la n(t_1-t_0)N^{3-2\al}\int_0^1(1-x)x^{1-2\al}\D x
=\frac{\la n^{H+1}(t_1-t_0)(x_1-x_0)^{2H}}{H(2H-1)}.
\end{align*}
Hence by Remark \ref{coefficient} we have, as $n\rightarrow\infty$, 
\begin{align*}
&{\mathbb E}_{{\mathbb Q}} [(n^{-\frac{H+1}{2}}{\mathcal S}_{1}^n(f))^2]\longrightarrow\frac{\la (t_1-t_0)(x_1-x_0)^{2H}}{H(2H-1)}\\
=&\la \int_0^1\int_{R^2}f(t,x)f(t,y)|x-y|^{2H-2}\D t\D x
={\mathbb E}_H\left(\int_0^1\int_Rf(t,x)W(\D t\D x)\right)^2.
\end{align*}
Meanwhile, substituting $\om({ i},{ x})=\sum_{y=-\infty}^\infty\psi_{y-x}\xi_{i,y}$ into ${\mathcal S}_{1}^n(f)$ yields
\begin{align}\label{lindeberg2}
{\mathcal S}_{1}^n(f)=2^{\frac{1}{2}}\sum_{y=-\infty}^\infty\sum_{i=1}^n\sum_{x\in{\cal E}_i}\psi_{y-x}\xi_{i,y},
\end{align}
for which one can use the method in Proposition \ref{clt1} to show $b'_{n,i,y}$ convergence to 0 uniformly. Hence in this case we have 
$$
n^{-\frac{H+1}{2}}{\mathcal S}_1^n(f)\stackrel{D}{\longrightarrow}\int_{[0,1]}\int_{R}
f({t},{x})W(\D{t}\D{x})=I_1(f).
$$

By Lemma \ref{chaosapproximation}, the remaining thing is to show Lemma \ref{k=1} true for $f$ of the form of linear combination of indicator functions $f_j,j=1,2,\dots,m,$ with disjoint,  finite area rectangle in $[0,1]\times \RR$. Actually, ${\mathcal S}_1^n(f_j), j=1,2,\dots,m$ are mutual independent when their time intervals are disjoint. Therefore we only consider the case of $m=2$ and the time interval  overlapping and space interval being disjoint for two disjoint rectangles. To this end, we assume 
\begin{align}
\label{indicator1}
f_1(t,x)=\mathbbm{1}_{[t_1, t_3]\times[x_1, x_2]}(t,x) ~{\text {and~}} f_2(t,x)=\mathbbm{1}_{[t_2, t_4]\times[x_3, x_4]}(t,x) 
\end{align}
for some $0\leq t_1\leq t_2\leq t_3\leq t_4\leq 1,x_1\leq x_2\leq x_3\leq x_4$.  Then, for $\la_1,\la_2\in{\mathbb R}$,
\begin{align*}
&{\mathcal S}_1^n(\la_1f_1+\la_2f_2)\\
=&2^{1/2}\la_1\!\!\!\sum_{\{nt_1\leq { i}\leq nt_3\}}\!\sum_{\{\sqrt{n}x_1\leq { x}\leq\sqrt{n}x_2\}}\!\!\!
\om({i},{x})\mathbbm{1}_{\{{ i}\leftrightarrow{ x}\}}+2^{1/2}\la_2\!\!\!\sum_{\{nt_2\leq { i}\leq nt_4\}}\sum_{\{\sqrt{n}x_3\leq {x}\leq\sqrt{n}x_4\}}\!\!\!
\om({i},{x})\mathbbm{1}_{\{{ i}\leftrightarrow{ x}\}}
\end{align*}
by linearity. We split ${\mathcal S}_1^n(\la_1f_1+\la_2f_2)$ into four terms:
\begin{align*}
&~~~~2^{1/2}\la_1\sum_{\{nt_1\leq { i}\leq nt_2\}}\sum_{\{\sqrt{n}x_1\leq { x}\leq\sqrt{n}x_2\}}
\om({i},{x})\mathbbm{1}_{\{{ i}\leftrightarrow{ x}\}}\\
&+2^{1/2}\la_1\sum_{\{nt_2< { i}< nt_3\}}\sum_{\{\sqrt{n}x_1\leq {x}\leq\sqrt{n}x_2\}}
\om({i},{x})\mathbbm{1}_{\{{ i}\leftrightarrow{ x}\}}\\
&+2^{1/2}\la_2\sum_{\{ nt_2< { i}< nt_3\}}\sum_{\{\sqrt{n}x_3\leq {x}\leq\sqrt{n}x_4\}}
\om({i},{x})\mathbbm{1}_{\{{ i}\leftrightarrow{ x}\}}\\
&+2^{1/2}\la_2\sum_{\{nt_3\leq { i}\leq nt_4\}}\sum_{\{\sqrt{n}x_3\leq {x}\leq\sqrt{n}x_4\}}
\om({i},{x})\mathbbm{1}_{\{{ i}\leftrightarrow{ x}\}}\\
=&J_1+J_{21}+J_{22}+J_3.
\end{align*}
By previous proof we know
$$
n^{-\frac{H+1}{2}}J_1\stackrel{D}{\longrightarrow}\la_1\int_{[0,1]}\int_{\RR}
\mathbbm{1}_{[t_1,t_2]\times[x_1,x_2]}(t,x)W(\D{t}\D{x}),
$$
$$
n^{-\frac{H+1}{2}}J_3\stackrel{D}{\longrightarrow}\la_2\int_{[0,1]}\int_{\RR}
\mathbbm{1}_{[t_3,t_4]\times[x_3,x_4]}(t,x)W(\D{t}\D{x})
$$
as $n\longrightarrow\infty$. Notice that ${\mathcal S}_1^n(\la_1f_1+\la_2f_2)$ is a independent sum of three terms $J_1,J_{21}+J_{22}$ and $J_3$ since the time interval is non-overlap. Put $\tilde{f}=\la_1\mathbbm{1}_{[t_2,t_3]\times[x_1,x_2]}+\la_2\mathbbm{1}_{[t_2,t_3]\times[x_3,x_4]}$. 
If we can show 
\begin{equation}\label{inter}
n^{-\frac{H+1}{2}}(J_{21}+J_{22})\stackrel{D}{\longrightarrow}\int_0^1\int_R\tilde{f}(t,x)W(\D t\D x)
\end{equation}
as $n\to\infty$, then we have
\begin{align*}
&n^{-\frac{H+1}{2}}{\mathcal S}_1^n(\la_1f_1+\la_2f_2)\\
=&n^{-\frac{H+1}{2}}J_1+n^{-\frac{H+1}{2}}(J_{21}+J_{22})+n^{-\frac{H+1}{2}}J_3\\
&\stackrel{D}{\longrightarrow}\int_0^1\int_{\RR}(\la_1f_1+\la_2f_2)(t,x)W(\D t\D x).
\end{align*}
 We compute the expectation
\begin{align*}
&{\mathbb E}_{{\mathbb Q}}(J_{21}J_{22})\\
=&2\la_1\la_2{\mathbb E}_{{\mathbb Q}}\bigg(\sum_{nt_2< { i}< nt_3}\sum_{\sqrt{n}x_1\leq {x}\leq\sqrt{n}x_2}
\sum_{ nt_2< {j}< nt_3}\sum_{\sqrt{n}x_3\leq {y}\leq\sqrt{n}x_4}
\om({i},{x})\om({j},{y})\mathbbm{1}_{\{{ i}\leftrightarrow{ x}\}}1_{\{{ j}\leftrightarrow{ y}\}}\bigg)\\
=&2\la_1\la_2\sum_{nt_2< { i}< nt_3}\sum_{\sqrt{n}x_1\leq {x}\leq\sqrt{n}x_2}
\sum_{\sqrt{n}x_3\leq {y}\leq\sqrt{n}x_4}\ga(y-x)
\mathbbm{1}_{\{{ i}\leftrightarrow{ x}\}}1_{\{{ i}\leftrightarrow{ y}\}}\\
\sim&\la\la_1\la_2\sum_{nt_2< { i}< nt_3}\sum_{\sqrt{n}x_1\leq {x}\leq\sqrt{n}x_2}
\sum_{\sqrt{n}x_3\leq {y}\leq\sqrt{n}x_4}\left(\frac{y-x}{\sqrt{n}}\right)^{1-2\al}n^{\frac{3}{2}-\al}\frac{1}{n}\quad\quad (\mbox{Noticing~} x\leq y)\\
\sim&\la\la_1\la_2n^{\frac{5}{2}-\al}(t_3-t_2)\int_{x_1}^{x_2}\!\!\int_{x_3}^{x_4}|x-y|^{2H-2}\D x\D y\\
=&\la_1\la_2(t_3-t_2)\left[((x_4-x_1)^{2H}-(x_4-x_2)^{2H}+(x_3-x_2)^{2H}-(x_3-x_1)^{2H})\right]n^{1+H}.
\end{align*}
Also one has
\begin{align*}
|\tilde{f}|^2_{H}=&H(2H-1)\int_0^1\int_{\RR^2}\tilde{f}(t,x)\tilde{f}(t,y)|x-y|^{2H-2}\D t\D x\D y\\
=&\la_1^2H(2H-1)\int_{t_2}^{t_3}\!\!\int_{x_1}^{x_2}\!\!\int_{x_1}^{x_2}|x-y|^{2H-2}\D t\D x\D y\\
&+\la_2^2H(2H-1)\int_{t_2}^{t_3}\!\!\int_{x_3}^{x_4}\!\!\int_{x_3}^{x_4}|x-y|^{2H-2}\D t\D x\D y\\
&+2\la_1\la_2H(2H-1)\int_{t_2}^{t_3}\!\!\int_{x_1}^{x_2}\!\!\int_{x_3}^{x_4}|x-y|^{2H-2}\D t\D x\D y\\
=&(t_3-t_2)\left[\la_1^2(x_2-x_1)^{2H}+\la_2^2(x_4-x_3)^{2H}\right]\\
&+2\la_1\la_2(t_3-t_2)\left[(x_4-x_1)^{2H}-(x_4-x_2)^{2H}+(x_3-x_2)^{2H}-(x_3-x_1)^{2H}\right].
\end{align*}
It follows that 
$$n^{-(H+1)}{\mathbb E}_{{\mathbb Q}}[(J_{21}+J_{22})^2]\longrightarrow |\tilde{f}|^2_{H}={\mathbb E}_H\left[\int_0^1\!\int_{\RR}\tilde{f}(t,x)W(\D t\D x)\right]^2\quad\quad \mbox{as~} n\longrightarrow\infty.
$$
For the sake of verifying Lindeberg's condition, we rearrange $J_{21}+J_{22}$ and have 
\begin{align*}
J_{21}+J_{22}=\sum_{y=-\infty}^\infty2^{1/2}\sum_{nt_2<i<nt_3}\bigg(\sum_{\stackrel{\sqrt{n}x_1\leq x\leq\sqrt{n}x_2}{i\leftrightarrow x}}\la_1\psi_{y-x}+\!\!\sum_{\stackrel{\sqrt{n}x_3\leq x\leq\sqrt{n}x_4}{i\leftrightarrow x}}\la_2\psi_{y-x}\bigg)\xi_{i,y},
\end{align*}
and 
\begin{align*}
{\mathbb E}_{{\mathbb Q}}[(J_{21}+J_{22})^2]&=2\sum_{y=-\infty}^\infty\sum_{nt_2<i<nt_3}\bigg(\sum_{\stackrel{\sqrt{n}x_1\leq x\leq\sqrt{n}x_2}{i\leftrightarrow x}}\la_1\psi_{y-x}+\!\!\sum_{\stackrel{\sqrt{n}x_3\leq x\leq\sqrt{n}x_4}{i\leftrightarrow x}}\la_2\psi_{y-x}\bigg)^2\\
&:=2\sum_{y=-\infty}^\infty\sum_{nt_2<i<nt_3}|\Psi(n,i,y)|^2.
\end{align*}
Since $[x_1,x_2]\cap[x_3,x_4]$ is empty, one has, by Cauchy-Schwarz inequality,
\begin{align*}
&\frac{1}{{\mathbb E}_{{\mathbb Q}}[(J_{21}+J_{22})^2]}|\Psi(n,i,y)|^2\\
\leq&\frac{n^{1/2}[\la_1^2(x_2-x_1)+\la_2^2(x_4-x_3)]}{{\mathbb E}_{{\mathbb Q}}[(J_{21}+J_{22})^2]}\sum_{\stackrel{\sqrt{n}x_1\leq x\leq\sqrt{n}x_2}{\sqrt{n}x_3\leq x\leq\sqrt{n}x_4}}\psi^2_{y-x}\\
\leq&\frac{n^{1/2}[\la_1^2(x_2-x_1)+\la_2^2(x_4-x_3)]}{{\mathbb E}_{{\mathbb Q}}[(J_{21}+J_{22})^2]}\sum_{x=-\infty}^\infty\psi^2_{x},
\end{align*}
which tends to zero uniformly in $y$ as $n\longrightarrow\infty$. Hence we can proceed as the proof in proposition \ref{clt1} to show \eqref{inter} holding. Furthermore, by a density argument we can confirm that \eqref{mainthm} is true for the case of $k=1$.

In addition, the above proof also implies that  
$$(n^{-\frac{H+1}{2}}{\mathcal S}_1^n(f_1),n^{-\frac{H+1}{2}}{\mathcal S}_1^n(f_2),\dots,n^{-\frac{H+1}{2}}{\mathcal S}_1^n(f_m))
\stackrel{D}{\longrightarrow}(I_1(f_1),I_1(f_2),\dots,I_1(f_m))
$$
holds as $n\longrightarrow\infty$ by Cram{\'e}r-Wold device and density arguments for $f_1,f_2,\dots,f_m\in {\mathcal L}_H$.\qed

We proceed to  go through with the case of $k>1$. Since no such technique as deletion of diagonals due to the spatial-colored property of fractional Brownian motion is  available we are not allowed to use the same techniques as in \cite{alberts2014}.
\vskip 0.1in 

{\it Proof of Lemma \ref{k>1}}.\quad We need to show
\begin{align}\label{k-multiple}
n^{-\frac{(H+1)k}{2}}{\mathcal S}_k^n(g^{\otimes k})\longrightarrow I_k(g^{\otimes k})\quad {\text as}\quad n\longrightarrow\infty.
\end{align}

We  begin with $k=2$. In this case, noticing that no same time indices appear in every term in ${\mathcal S}_2^n$, we have that ${\mathcal S}_2^n(g^{\otimes 2})$ equals
\begin{align}\label{2discrete}
&2\sum_{{\mathbf i}\in E_2^n}\sum_{{\mathbf x}\in {\mathbb Z}^2}
\mathbbm{1}^{\otimes 2}_{[t_0,t_1]\times [x_0,x_1]}\Big(\frac{{\mathbf i}}{n},\frac{{\mathbf x}}{\sqrt{n}}\Big)\om({\mathbf i},{\mathbf x})\mathbbm{1}_{\{{\mathbf i}\leftrightarrow{\mathbf x}\}}\nonumber\\
=&[{\mathcal S}_1^n({g})]^2-2\sum_{{ i}\in E_1^n}\sum_{{\mathbf x}\in {\mathbb Z}^2}
\mathbbm{1}_{[t_0, t_1]}\Big(\frac{{i}}{n}\Big)\mathbbm{1}^{\otimes 2}_{[x_0,x_1]}\Big(\frac{{\mathbf x}}{\sqrt{n}}\Big)\om({i},{x_1})\om({i},{x_2})\mathbbm{1}_{\{{i}\leftrightarrow{x_1}\}}\mathbbm{1}_{\{{i}\leftrightarrow{x_2}\}}\nonumber\\
=&[{\mathcal S}_1^n(g)]^2-2\sum_{nt_0\leq { i}\leq nt_1}\Big(\sum_{{ x}\in {\mathbb Z}}
\mathbbm{1}_{[x_0,x_1]}\Big(\frac{{ x}}{\sqrt{n}}\Big)\om({i},{x})\mathbbm{1}_{\{{i}\leftrightarrow{x}\}}\Big)^2.
\end{align}
By continuous mapping theorem, we know 
\begin{align}\label{conmap}
n^{-(H+1)}[{\mathcal S}_1^n({g})]^2\stackrel{D}{\longrightarrow} \left(\int_0^1\int_{\mathbb R}
\mathbbm{1}_{[t_0,t_1]\times[x_0,x_1]}(t,x)W(\D t\D x)\right)^2.
\end{align}
Put $J_i:=\sum_{{ x}\in {\mathbb Z}}
\mathbbm{1}_{[x_0,x_1]}\big(\frac{{ x}}{\sqrt{n}}\big)\om({i},{x})\mathbbm{1}_{\{{i}\leftrightarrow{x}\}}, i=\lfloor nt_0\rfloor+1,\dots,\lfloor nt_1\rfloor$, then $J_i$'s are independent. Since the variance of $\sqrt{2}J_i$ is 
$$
N\ga(0)+2\sum_{k=1}^{N-1}(N-k)\ga(k)\sim N\ga(0)+2N^{3-2\al}\la\sum_{k=1}^{N-1}\left(1-\frac{k}{N}\right)\left(\frac{k}{N}\right)^{1-2\al}\frac{1}{N}\sim N\ga(0)+N^{3-2\al}.
$$
 As before, we can show
$$
\sqrt{2}n^{-\frac{H}{2}}J_i\stackrel{D}{\longrightarrow}N(0,(x_1-x_0)^{2H})
$$
as $n\longrightarrow\infty$.

Therefore, by the law of large number for triangular arrays (see \cite[Theorem 2.3.4]{Durrett1996Probability}), we have
\begin{align}\label{2ergodicity}
&~~~~2n^{-(H+1)}\sum_{nt_0\leq { i}\leq nt_1}\bigg(\sum_{{ x}\in {\mathbb Z}}
\mathbbm{1}_{[x_0,x_1]}\left(\frac{{ x}}{\sqrt{n}}\bigg)\om({i},{x})\mathbbm{1}_{\{{i}\leftrightarrow{x}\}}\right)^2\nonumber\\
&=2n^{-(H+1)}\sum_{nt_0\leq { i}\leq nt_1}\bigg(\sum_{{ x}\in {\mathbb Z}}
\mathbbm{1}_{[x_0,x_1]}\left(\frac{{ x}}{\sqrt{n}}\right)\om({i},{x})\mathbbm{1}_{\{{i}\leftrightarrow{x}\}}\bigg)^2\nonumber\\
&\stackrel{P}{\longrightarrow}(t_1-t_0)(x_1-x_0)^{2H}\quad\text{as}~n\longrightarrow\infty.
\end{align}

Put $\tilde{g}=g/(t_1-t_0)(x_1-x_0)^{2H}$, then $\|\tilde{g}\|_H=1$, \eqref{2discrete}, \eqref{conmap}, \eqref{2ergodicity} and \eqref{contract2} imply
$$
n^{-(H+1)}{\mathcal S}_2^n(\tilde{g}^{\otimes 2})\stackrel{D}{\longrightarrow}\left(\int_0^1\int_{\mathbb R}
\tilde{g}(t,x)W(\D t\D x)\right)^2-1=H_2(W(\tilde{g}))=I_2(\tilde{g}^{\otimes 2}).
$$

Now assume that \eqref{k-multiple} holds for $k=1,2,\dots,l-1$. Rewrite ${\mathcal S}_l^n({f})$, with $g(t,x)=\mathbbm{1}_{[t_0,t_1]\times [x_0,x_1]}(t,x)$, as 
\begin{equation*}
{\mathcal S}_l^n(g^{\otimes l})
=2^{l/2}\sum_{{\mathbf i}\in E_l^n}\sum_{{\mathbf x}\in {\mathbb Z}^l}
\overline {g_n^{\otimes l}}\left(\frac{{\mathbf i}}{n},\frac{{\mathbf x}}{\sqrt{n}}\right)\om({\mathbf i},{\mathbf x})\mathbbm{1}_{\{{\mathbf i}\leftrightarrow{\mathbf x}\}},
\end{equation*}
which also can be expressed as the difference of 
\begin{align}\label{S1}
2^{\frac{l}{2}}&\bigg(\sum_{{\bar{\mathbf i}}\in E_{l-1}^n}\sum_{{\mathbf y}\in {\mathbb Z}^{l-1}}
\overline {g_n^{\otimes (l-1)}}\left(\frac{{\mathbf {\bar i}}}{n},\frac{{\mathbf y}}{\sqrt{n}}\right)\om({\mathbf {\bar i}},{\mathbf y})\mathbbm{1}_{\{{\mathbf {\bar i}}\leftrightarrow{\mathbf y}\}}\bigg)\nonumber\\
&~~~~~\times\bigg(\sum_{{i}_l\in E_1^n}\sum_{{x}_l\in {\mathbb Z}}
{\bar g}_n\left(\frac{{ i}_l}{n},\frac{{x}_l}{\sqrt{n}}\right)\om({i}_l,{ x_l})\mathbbm{1}_{\{{i_l}\leftrightarrow{x_l}\}}\bigg)
\end{align}
and $R_{l,n}$. Here ${\mathbf {\bar i}}=(i_1,\dots,i_{l-1})$, $R_{l,n}$ is the sum of all terms that they are multiplication of two factors, which comes from the two parentheses of \eqref{S1}, respectively, with exactly one component in ${\mathbf {\bar i}}$ corresponding to the first factor equal to $i_l$ corresponding to the second factor. Specifically, denote by ${\mathbf {\bar i}}^j$ the vector by canceling the jth component $i_j$ of ${\mathbf {\bar i}}$, and observe that all components in ${\mathbf {\bar i}}^j$ are different from $i_j$, $j=1,2,\dots l-1$. It is necessary to remark that the dimension of vector ${\mathbf i}$ may vary in the different display, say, the second line of \eqref{S3}. Then $2^{-\frac{l}{2}}R_{l,n}$ equals
\begin{equation*}
\begin{split}
&\sum_{i_1\in A}\bigg[\sum_{x\in\ZZ}{\bar g}_n\left(\frac{{i_1}}{n},\frac{{ x}}{\sqrt{n}}\right)\om({i_1},{x})\mathbbm{1}_{\{{i_1}\leftrightarrow{x}\}}\bigg]^2\sum_{{\mathbf {\bar i}}^1\in E_{l-2}^n}\sum_{{\mathbf x}\in {\mathbb Z}^{l-2}}
\overline { g^{\otimes (l-2)}_n}\left(\frac{\bar{{\mathbf i}}^1}{n},\frac{{\mathbf x}}{\sqrt{n}}\right)\om(\bar{{\mathbf i}}^1,{\mathbf x})\mathbbm{1}_{\{\bar{{\mathbf i}}^1\leftrightarrow{\mathbf x}\}}\\
+&\sum_{i_2\in A}\bigg[\sum_{x\in\ZZ}{\bar g}_n\left(\frac{{i_2}}{n},\frac{{ x}}{\sqrt{n}}\right)\om({i_2},{x})\mathbbm{1}_{\{{i_2}\leftrightarrow{x}\}}\bigg]^2\sum_{{\mathbf {\bar i}}^2\in E_{l-2}^n}\sum_{{\mathbf x}\in {\mathbb Z}^{l-2}}
\overline { g^{\otimes (l-2)}_n}\left(\frac{\bar{{\mathbf i}}^2}{n},\frac{{\mathbf x}}{\sqrt{n}}\right)\om(\bar{{\mathbf i}}^2,{\mathbf x})\mathbbm{1}_{\{\bar{{\mathbf i}}^2\leftrightarrow{\mathbf x}\}}\\
+&\cdots
\end{split}
\end{equation*}
\begin{equation}\label{S2}
\begin{split}
+&\sum_{i_{l-1}\in A}\bigg(\bigg[\sum_{x\in\ZZ}{\bar g}_n\left(\frac{{i_{l-1}}}{n},\frac{{ x}}{\sqrt{n}}\right)\om({i_{l-1}},{x})\mathbbm{1}_{\{{i_{l-1}}\leftrightarrow{x}\}}\bigg]^2\\
&\quad\quad\quad\quad\quad\quad\quad\quad\quad\quad\quad\quad\quad\times\sum_{{\mathbf {\bar i}}^{l-1}\in E_{l-2}^n}\sum_{{\mathbf x}\in {\mathbb Z}^{l-2}}
\overline { g^{\otimes (l-2)}_n}\left(\frac{\bar{{\mathbf i}}^{l-1}}{n},\frac{{\mathbf x}}{\sqrt{n}}\right)\om(\bar{{\mathbf i}}^{l-1},{\mathbf x})\mathbbm{1}_{\{\bar{{\mathbf i}}^{l-1}\leftrightarrow{\mathbf x}\}}\bigg)\\
&=2(l-1)\sum_{i_1\in A}\bigg(\bigg[\sum_{x\in\ZZ}{\bar g}_n\left(\frac{{i_1}}{n},\frac{{ x}}{\sqrt{n}}\right)\om({i_1},{x})\mathbbm{1}_{\{{i_1}\leftrightarrow{x}\}}\bigg]^2\\
&\quad\quad\quad\quad\quad\quad\quad\quad\quad\quad\quad\quad\quad\times\sum_{{\mathbf {\bar i}}^1\in E_{l-2}^n}\sum_{{\mathbf x}\in {\mathbb Z}^{l-2}}
\overline {g^{\otimes(l-2)}_n}\left(\frac{\bar{{\mathbf i}}^1}{n},\frac{{\mathbf x}}{\sqrt{n}}\right)\om(\bar{{\mathbf i}}^1,{\mathbf x})\mathbbm{1}_{\{\bar{{\mathbf i}}^1\leftrightarrow{\mathbf x}\}}\bigg),
\end{split}
\end{equation}
where the set $A$ consists of all integers in the interval $[nt_0,nt_1]$.
 Put 
\begin{equation}\label{approx}
{\mathcal S}_{l-2}^{n,1}(g^{\otimes (l-2)}):=2^{\frac{l-2}{2}}\sum_{{\mathbf {\bar i}}^1\in E_{l-2}^n}\sum_{{\mathbf x}\in {\mathbb Z}^{l-2}}
\overline { g^{\otimes (l-2)}_n}\bigg(\frac{\bar{{\mathbf i}}^1}{n},\frac{{\mathbf x}}{\sqrt{n}}\bigg)\om(\bar{{\mathbf i}}^1,{\mathbf x})\mathbbm{1}_{\{\bar{{\mathbf i}}^1\leftrightarrow{\mathbf x}\}}
\end{equation}
with $i_1\in [n]$ fixed. We are going to estimate the difference between ${\mathcal S}_{l-2}^{n,1}(g^{\otimes (l-2)})$ and ${\mathcal S}_{l-2}^{n}(g^{\otimes (l-2)})$. Since the time index ${\mathbf i}\in [n]^{k-2}$ of the summand which is in the sum of ${\mathcal S}_{l-2}^{n}(g^{\otimes (l-2)})$
but not in that of ${\mathcal S}_{l-2}^{n,1}(g^{\otimes (l-2)})$ must have one and only one component equal to $i_1$, we have
\begin{equation}\label{S3}
\begin{split}
&{\mathcal S}_{l-2}^{n}(g^{\otimes (l-2)})-{\mathcal S}_{l-2}^{n,1}(g^{\otimes (l-2)})\\
=&(l-2)2^{\frac{l-2}{2}}\sum_{y\in \ZZ}{\bar g}_n\bigg(\frac{{i_1}}{n},\frac{{ y}}{\sqrt{n}}\bigg)\om({i_1},{y})\mathbbm{1}_{\{{i_1}\leftrightarrow{y}\}}\sum_{{\mathbf {\bar i}}^1\in E_{l-3}^n}\sum_{{\mathbf x}\in {\mathbb Z}^{l-3}}
\overline { g^{\otimes (l-3)}_n}\bigg(\frac{\bar{{\mathbf i}}^1}{n},\frac{{\mathbf x}}{\sqrt{n}}\bigg)\om(\bar{{\mathbf i}}^1,{\mathbf x})\mathbbm{1}_{\{\bar{{\mathbf i}}^1\leftrightarrow{\mathbf x}\}}.
\end{split}
\end{equation}
Hence we get from \eqref{S2} and \eqref{S3} that $R_{l,n}$ is the difference of 
$$
2(l-1){\mathcal S}_{l-2}^{n}(g^{\otimes (l-2)})\sum_{i_1\in A}\bigg[\sum_{x\in\ZZ}{\bar g}_n\bigg(\frac{{i_1}}{n},\frac{{ x}}{\sqrt{n}}\bigg)\om({i_1},{x})\mathbbm{1}_{\{{i_1}\leftrightarrow{x}\}}\bigg]^2
$$
and 
\begin{align*}
(l-1)(l-2)2^{\frac{l}{2}}\sum_{i_1\in A}&\bigg[\sum_{x\in\ZZ}{\bar g}_n\bigg(\frac{{i_1}}{n},\frac{{ x}}{\sqrt{n}}\bigg)\om({i_1},{x})\mathbbm{1}_{\{{i_1}\leftrightarrow{x}\}}\bigg]^3\\
&\times\sum_{{\mathbf {\bar i}}^1\in E_{l-3}^n}\sum_{{\mathbf x}\in {\mathbb Z}^{l-3}}
{\bar g}_n\bigg(\frac{\bar{{\mathbf i}}^1}{n},\frac{{\mathbf x}}{\sqrt{n}}\bigg)\om(\bar{{\mathbf i}}^1,{\mathbf x})\mathbbm{1}_{\{\bar{{\mathbf i}}^1\leftrightarrow{\mathbf x}\}}.
\end{align*}
Let $Q_{i_1}=\sum_{x\in\ZZ}{\bar g}_n\big(\frac{{i_1}}{n},\frac{{ x}}{\sqrt{n}}\big)\om({i_1},{x})\mathbbm{1}_{\{{i_1}\leftrightarrow{x}\}}$ for $i_1\in A$, then $R_{l,n}$ can be written as
\begin{align*}
2(l-1){\mathcal S}_{l-2}^{n}&(g^{\otimes (l-2)})\sum_{i_1\in A}Q_{i_1}^2\\
-&(l-1)(l-2)2^{\frac{l}{2}}\sum_{i_1\in A}Q_{i_1}^3\sum_{{\mathbf {\bar i}}^1\in E_{l-3}^n}\sum_{{\mathbf x}\in {\mathbb Z}^{l-3}}
{\bar g}_n\bigg(\frac{\bar{{\mathbf i}}^1}{n},\frac{{\mathbf x}}{\sqrt{n}}\bigg)\om(\bar{{\mathbf i}}^1,{\mathbf x})\mathbbm{1}_{\{\bar{{\mathbf i}}^1\leftrightarrow{\mathbf x}\}}.
\end{align*}

Considering the term $\sum\limits_{{\mathbf {\bar i}}^1\in E_{l\!\!-3}^n}\!\!\sum\limits_{{\mathbf x}\in {\mathbb Z}^{l\!-3}}\!\!
{\bar g}_n\big(\frac{\bar{{\mathbf i}}^1}{n},\frac{{\mathbf x}}{\sqrt{n}}\big)\om(\bar{{\mathbf i}}^1,{\mathbf x})\mathbbm{1}_{\{\bar{{\mathbf i}}^1\leftrightarrow{\mathbf x}\}}
$ and keeping the process as before, we obtain that $R_{l,n}$
\begin{equation}\label{S4}
\begin{split}
=&(l-1)2{\mathcal S}_{l-2}^{n}(g^{\otimes (l-2)})\sum_{i_1\in A}Q_{i_1}^2-2^{\frac{3}{2}}(l-1)(l-2){\mathcal S}_{l-3}^{n}(g^{\otimes (l-3)})\sum_{i_1\in A}Q_{i_1}^3\\
&-2^{\frac{l}{2}}(l-1)(l-2)(l-3)\sum_{i_1\in A}Q_{i_1}^4\sum_{{\mathbf {\bar i}}^1\in E_{l-4}^n}\sum_{{\mathbf x}\in {\mathbb Z}^{l-4}}
\overline {g_n^{\otimes(l-4)}}\bigg(\frac{\bar{{\mathbf i}}^1}{n},\frac{{\mathbf x}}{\sqrt{n}}\bigg)\om(\bar{{\mathbf i}}^1,{\mathbf x})\mathbbm{1}_{\{\bar{{\mathbf i}}^1\leftrightarrow{\mathbf x}\}}=\cdots\\
=&(l-1)2{\mathcal S}_{l-2}^{n}(g^{\otimes (l-2)})\sum_{i_1\in A}Q_{i_1}^2-2^{\frac{3}{2}}(l-1)(l-2){\mathcal S}_{l-3}^{n}(g^{\otimes (l-3)})\sum_{i_1\in A}Q_{i_1}^3\\
&+2^2(l-1)(l-2)(l-3){\mathcal S}_{l-4}^{n}(g^{\otimes (l-4)})\sum_{i_1\in A}Q_{i_1}^4+\cdots+(-\sqrt{2})^{l-1}(l-1)!{\mathcal S}_{1}^{n}(g)\sum_{i_1\in A}Q_{i_1}^{l-1}\\
&+(-\sqrt{2})^{l}(l-1)!\sum_{i_1\in A}Q_{i_1}^l.
\end{split}
\end{equation}
Multiplying $R_{l,n}$ with the coefficient $n^{-\frac{(H+1)l}{2}}$, we can find that the first term in the RHS of the last equal sign above converges to 
$(l-1)I_{l-2}(h^{\otimes (l-2)})|h|_H^2$ and all other terms converge to zero as $n\longrightarrow\infty$. Actually, as \eqref{2ergodicity}, we have 
$$
\frac{1}{\lfloor n(t_1-t_0)\rfloor}\sum_{i_1\in A}(Q_{i_1}^2-\EE_{\mathbb Q} Q_{i_1}^2)\longrightarrow 0
$$
as $n\longrightarrow \infty$. While
$2\EE_{{\mathbb Q}} Q_{i_1}^2\sim N\ga(0)+2\la\sum_{i=1}^N(N-i)i^{1-2\al}$
$\sim N\ga(0)+\frac{\la N^{2H}}{H(2H-1)}$. Hence we get 
$$
2n^{-(H+1)}\sum_{i_1\in A}Q_{i_1}^2\longrightarrow(t_1-t_0)(x_1-x_0)^{2H}, ~~\text{and}~~
$$
$$
2(l-1)n^{-\frac{(H+1)l}{2}}{\mathcal S}_{l-2}^{n}(g^{\otimes (l-2)})\sum_{i_1\in A}Q_{i_1}^2\stackrel{D}{\longrightarrow}(l-1)I_{l-2}(g^{\otimes (l-2)})(t_1-t_0)(x_1-x_0)^{2H},
$$
as $n\longrightarrow\infty.$  

Now we only check the third term in $R_{l,n}$ (the other is similar). Since 
$$
n^{-2}\sum_{i_1\in A}\bigg(\frac{Q_{i_1}}{n^{H/2}}\bigg)^4\longrightarrow 0
$$ 
as \eqref{2ergodicity} we have
\begin{align*}
&2^2(l-1)(l-2)(l-3)n^{-\frac{(H+1)l}{2}}{\mathcal S}_{l-4}^{n}(g^{\otimes (l-4)})\sum_{i_1\in A}Q_{i_1}^4\\
=&2^2(l-1)(l-2)(l-3)n^{-\frac{(H+1)(l-4)}{2}}{\mathcal S}_{l-4}^{n}(g^{\otimes (l-4)})n^{-2}\sum_{i_1\in A}\bigg(\frac{Q_{i_1}}{n^{H/2}}\bigg)^4\longrightarrow 0.
\end{align*}

Then, by induction hypothesis and using \eqref{contract2}, \eqref{S1} and \eqref{S4}, we arrive at 
\begin{align*}
&n^{-\frac{(H+1)l}{2}}{\mathcal S}_l^n(g^{\otimes l})\\
\longrightarrow&I_{l-1}(g^{\otimes(l-1)})I_1({g})-(l-1)I_{l-2}(g^{\otimes (l-2)})(t_1-t_0)(x_1-x_0)^{2H}=I_l(g^{\otimes l}),
\end{align*}
as $n\longrightarrow\infty$ since $g^{\otimes l}\otimes_1g=\|g\|_H^2g^{\otimes (l-2)}$ is obvious by the definition of contraction of two functions.\qed
\vskip 0.1in

{\it Proof of Lemma \ref{tensor}}.\quad We show it only for $s=2$ and $k_1=m,k_2=l$. We use $f,g$ instead of $g_1,g_2$ and assume $f, g$ be as \eqref{indicator}. Actually, we can reduce it to the case of  $f(t,x)=1_A(t)f_1(x), g(t,x)=1_A(t)g_1(x)$ for some interval $A\subset[0,1]$ and $f_1,g_1$ being indicators of two (disjoint) intervals in $\RR$ since $f, g$ have the same time supports. 
\begin{align*}
&n^{\frac{-(H+1)(m+l)}{2}}{\mathcal S}_{m+l}^n(f^{\otimes m}\otimes g^{\otimes l})\\
=&n^{\frac{-m(H+1)}{2}}{\mathcal S}_{m}^n(f^{\otimes m})n^{\frac{-l(H+1)}{2}}{\mathcal S}_{l}^n(g^{\otimes l})\stackrel{D}{\longrightarrow}I_m(f^{\otimes m})I_l(g^{\otimes l})=I_{m+l}(f^{\otimes m}\otimes g^{\otimes l})
\end{align*}
 is true, otherwise finer partition and linear property of symmetrical tensors  can work. 

We show it by induction. By Lemma \ref{k>1}, we know that it holds when $m=0$ or $l=0$. Hence, let $m$ be fixed and suppose 
\begin{equation}\label{tensor2}
\begin{split}
&n^{\frac{-(H+1)(r+s)}{2}}{\mathcal S}_k^n(f^{\otimes r}\otimes g^{\otimes s})\\
&\stackrel{D}{\longrightarrow}\int_{[0,1]^k}\int_{\RR^k}
f^{\otimes r}\otimes g^{\otimes s}({\mathbf t},{\mathbf x})W^{\otimes (r+s)}(\D{\mathbf t}\D{\mathbf x})=I_{r+s}(f^{\otimes r}\otimes g^{\otimes s})
\end{split}
\end{equation}
holding for all $r=1,2,\dots,m, s=1,2,\dots l-1$. By using a similar method, we have
\begin{equation}\label{S5}
\begin{split}
&2^{-\frac{m+l}{2}}{\mathcal S}_{m+l}^n(f^{\otimes m}\otimes g^{\otimes l})\\
=&\sum_{(\mathbf{ i,j})\in E_{m+l}^n}\sum_{(\mathbf{x,y})\in{\mathbb Z}^{m+l}}\overline {f_n^{\otimes m}}\bigg(\frac{{\mathbf i}}{n},\frac{{\mathbf x}}{\sqrt{n}}\bigg)\om({\mathbf i},{\mathbf x})\mathbbm{1}_{\{{\mathbf i}\leftrightarrow{\mathbf x}\}}\overline {g_n^{\otimes l}}\bigg(\frac{{\mathbf j}}{n},\frac{{\mathbf y}}{\sqrt{n}}\bigg)\om({\mathbf j},{\mathbf y})\mathbbm{1}_{\{{\mathbf j}\leftrightarrow{\mathbf y}\}}\\
=&\sum_{(\mathbf{ i,j})\in E_{m+l-1}^n}\sum_{(\mathbf{x,y})\in{\mathbb Z}^{m+l-1}}\overline {f_n^{\otimes m}}\bigg(\frac{{\mathbf i}}{n},\frac{{\mathbf x}}{\sqrt{n}}\bigg)\om({\mathbf i},{\mathbf x})\mathbbm{1}_{\{{\mathbf i}\leftrightarrow{\mathbf x}\}}\overline {g_n^{\otimes (l-1)}}\bigg(\frac{{\mathbf j}}{n},\frac{{\mathbf y}}{\sqrt{n}}\bigg)\om({\mathbf j},{\mathbf y})\mathbbm{1}_{\{{\mathbf j}\leftrightarrow{\mathbf y}\}}\\
&~~~~\times\sum_{k=1}^n\sum_{y\in{\mathbb Z}}{\bar g}_n\bigg(\frac{{k}}{n},\frac{{y}}{\sqrt{n}}\bigg)\om({k},{y})1_{\{k\leftrightarrow{ y}\}}-2^{-\frac{m+l}{2}}R_{(m+l),n}.
\end{split}
\end{equation}
A careful arrangement  yields that $R_{(m+l),n}$ is given by
\begin{align*}
&m2^{\frac{m+l}{2}}\sum_{k=1}^n\bigg(\sum_{y\in{\mathbb Z}}{\bar g}_n\bigg(\frac{{k}}{n},\frac{{y}}{\sqrt{n}}\bigg)\om({k},{y})\mathbbm{1}_{\{k\leftrightarrow{ y}\}}
\sum_{y\in{\mathbb Z}}{\bar f}_n\bigg(\frac{{k}}{n},\frac{{y}}{\sqrt{n}}\bigg)\om({k},{y})\mathbbm{1}_{\{k\leftrightarrow{ y}\}}\bigg)\\
&~~~~~~~~~~~~~~~~~~~~~~\times\sum_{(\bar{\mathbf{ i}}^k,\bar{\mathbf{j}}^k)\in E_{m+l-2}^n}\sum_{\{\mathbf{x}\in{\mathbb Z}^{m-1},{\mathbf y}\in{\mathbb Z}^{l-1}\}}\!\!\!\overline{f_n^{\otimes(m-1)}}\bigg(\frac{\bar{{\mathbf i}}^k}{n},\frac{{\mathbf x}}{\sqrt{n}}\bigg)\om(\bar{{\mathbf i}}^k,{\mathbf x})\mathbbm{1}_{\{{\bar{\mathbf i}}^k\leftrightarrow{\mathbf x}\}}\\
&~~~~~~~~~~~~~~~~~~~~~~\times\overline{g_n^{\otimes(l-1)}}\bigg(\frac{\bar{\mathbf{j}}^k}{n},\frac{{\mathbf y}}{\sqrt{n}}\bigg)\om(\bar{\mathbf{j}}^k,{\mathbf y})\mathbbm{1}_{\{\bar{\mathbf{j}}^k\leftrightarrow{\mathbf y}\}}\\
+&(l-1)2^{\frac{m+l}{2}}\sum_{k=1}^n\bigg(\sum_{y\in{\mathbb Z}}{\bar g}_n\bigg(\frac{{k}}{n},\frac{{y}}{\sqrt{n}}\bigg)\om({k},{y})\mathbbm{1}_{\{k\leftrightarrow{ y}\}}\bigg)^2\\
\times&\sum_{(\bar{\mathbf{ i}}^k,\bar{\mathbf{j}}^k)\in E_{m+l-2}^n}\sum_{\{\mathbf{x}\in{\mathbb Z}^{m},{\mathbf y}\in{\mathbb Z}^{l-2}\}}\overline{f_n^{\otimes m}}\bigg(\frac{\bar{\mathbf{ i}}^k}{n},\frac{{\mathbf x}}{\sqrt{n}}\bigg)\om(\bar{\mathbf{ i}}^k,{\mathbf x})\mathbbm{1}_{\{\bar{\mathbf{ i}}^k\leftrightarrow{\mathbf x}\}}\overline{g_n^{\otimes(l-2)}}\bigg(\frac{\bar{\mathbf{j}}^k}{n},\frac{{\mathbf y}}{\sqrt{n}}\bigg)\om(\bar{\mathbf{j}}^k,{\mathbf y})\mathbbm{1}_{\{\bar{\mathbf{j}}^k\leftrightarrow{\mathbf y}\}}.
\end{align*}
Here for the sake of easy reading we write $\{\mathbf{x}\in{\mathbb Z}^{m-1},{\mathbf y}\in{\mathbb Z}^{l-1}\}$ instead of $({\mathbf x,\mathbf y})\in {\mathbb Z}^{m+l-2}$,  
hence one can find quickly the dimension of $\bar{\mathbf i}, \bar{\mathbf j}$. Let 
$$Q_{n,k}=\sum_{y\in{\mathbb Z}}{\bar f}_n\bigg(\frac{{k}}{n},\frac{{y}}{\sqrt{n}}\bigg)\om({k},{y})\mathbbm{1}_{\{k\leftrightarrow{ y}\}}, P_{n,k}=\sum_{y\in{\mathbb Z}}{\bar g}_n\bigg(\frac{{k}}{n},\frac{{y}}{\sqrt{n}}\bigg)\om({k},{y})\mathbbm{1}_{\{k\leftrightarrow{ y}\}}.
$$
Since, for $\forall n,k\in\NN$, $\EE_{\QQ}(n^{-\frac{H}{2}}Q_{n,k})^2\leq C,\EE_{\QQ}(n^{-\frac{H}{2}}P_{n,k})^2\leq C$ for some $C>0$, then an application of Burkholder's inequality implies, for $\forall n,k,q\in\NN$, 
\begin{align}\label{boundedness}
\EE_{\QQ}(n^{-\frac{H}{2}}Q_{n,k})^q\leq C_q~~~\mbox{and}~~~\EE_{\QQ}(n^{-\frac{H}{2}}P_{n,k})^q\leq C_q
\end{align}
for some constant $C_q>0$.

A similar manner as in  \eqref{S4} shows that $R_{(m+l),n}$ equals the sum $B_1-B_2+B_3-B_4$. Here $B_1,B_3$ equal 
$$
2m{\mathcal S}_{m+l-2}^n(f^{\otimes (m-1)}\otimes g^{\otimes (l-1)})\sum_{k=1}^n(Q_{k,n}P_{k,n})
$$
and
$$
2(l-1){\mathcal S}_{m+l-2}^n(f^{\otimes m}\otimes g^{\otimes (l-2)})\sum_{k=1}^nP_{k,n}^2,
$$
respectively. $B_2$ equals $m2^{\frac{m+l}{2}}\sum_{k=1}^n(Q_{k,n}P_{k,n})$ times the sum of 
$$
\begin{aligned}
(m-1)Q_{k,n}&\sum_{(\bar{\mathbf{ i}}^k,{\bar{\mathbf{j}}^k})\in E_{m+l-3}^n}\sum_{(\mathbf{x,y})\in{\mathbb Z}^{m+l-3}}\overline{ f_n^{\otimes (m-2)}}\bigg(\frac{\bar{{\mathbf i}}^k}{n},\frac{{\mathbf x}}{\sqrt{n}}\bigg)\\
&\times\om(\bar{{\mathbf i}}^k,{\mathbf x})\mathbbm{1}_{\{{\bar{\mathbf i}}^k\leftrightarrow{\mathbf x}\}}\overline{ g_n^{\otimes (l-1)}}\bigg(\frac{\bar{\mathbf{j}}^k}{n},\frac{{\mathbf y}}{\sqrt{n}}\bigg)\om(\bar{\mathbf{j}}^k,{\mathbf y})\mathbbm{1}_{\{{\mathbf j}\leftrightarrow{\mathbf y}\}}
\end{aligned}
$$
and
$$
\begin{aligned}
(l-1)P_{k,n}&\sum_{(\bar{\mathbf{ i}}^k,\bar{\mathbf{j}}^k)\in E_{m+l-3}^n}\sum_{(\mathbf{x,y})\in{\mathbb Z}^{m+l-3}}\overline{ f_n^{\otimes (m-1)}}\bigg(\frac{\bar{{\mathbf i}}^k}{n},\frac{{\mathbf x}}{\sqrt{n}}\bigg)\\
&\times\om(\bar{{\mathbf i}}^k,{\mathbf x})\mathbbm{1}_{\{{\bar{\mathbf i}}^k\leftrightarrow{\mathbf x}\}}\overline{ g_n^{\otimes (l-2)}}\bigg(\frac{\bar{\mathbf{j}}^k}{n},\frac{{\mathbf y}}{\sqrt{n}}\bigg)\om(\bar{\mathbf{j}}^k,{\mathbf y})\mathbbm{1}_{\{\bar{\mathbf{j}}^k\leftrightarrow{\mathbf y}\}}.
\end{aligned}
$$
$B_4$ is given by the product of $(l-1)2^{\frac{m+l}{2}}\sum_{k=1}^nP_{k,n}^2$ and the sum of 
$$
\begin{aligned}
mQ_{k,n}&\sum_{(\bar{\mathbf{ i}}^k,{\bar{\mathbf{j}}^k})\in E_{m+l-3}^n}\sum_{(\mathbf{x,y})\in{\mathbb Z}^{m+l-3}}\overline{ f_n^{\otimes (m-1)}}\bigg(\frac{\bar{{\mathbf i}}^k}{n},\frac{{\mathbf x}}{\sqrt{n}}\bigg)\\
&\times\om(\bar{{\mathbf i}}^k,{\mathbf x})\mathbbm{1}_{\{{\bar{\mathbf i}}^k\leftrightarrow{\mathbf x}\}}\overline{ g_n^{\otimes (l-2)}}\bigg(\frac{\bar{\mathbf{j}}^k}{n},\frac{{\mathbf y}}{\sqrt{n}}\bigg)\om(\bar{\mathbf{j}}^k,{\mathbf y})\mathbbm{1}_{\{{\mathbf j}\leftrightarrow{\mathbf y}\}}
\end{aligned}
$$
and
$$
\begin{aligned}
(l-2)P_{k,n}&\sum_{(\bar{\mathbf{ i}}^k,\bar{\mathbf{j}}^k)\in E_{m+l-3}^n}\sum_{(\mathbf{x,y})\in{\mathbb Z}^{m+l-3}}\overline{ f_n^{\otimes m}}\bigg(\frac{\bar{{\mathbf i}}^k}{n},\frac{{\mathbf x}}{\sqrt{n}}\bigg)\\
&\times\om(\bar{{\mathbf i}}^k,{\mathbf x})\mathbbm{1}_{\{{\bar{\mathbf i}}^k\leftrightarrow{\mathbf x}\}}\overline{ f_n^{\otimes (l-3)}}\bigg(\frac{\bar{\mathbf{j}}^k}{n},\frac{{\mathbf y}}{\sqrt{n}}\bigg)\om(\bar{\mathbf{j}}^k,{\mathbf y})\mathbbm{1}_{\{\bar{\mathbf{j}}^k\leftrightarrow{\mathbf y}\}}.
\end{aligned}
$$
Repeating the procedure above to the sum in $B_2,B_4$, it follows 
\begin{equation}\label{S6}
\begin{split}
R_{(m+l),n}=&2m{\mathcal S}_{m+l-2}^n(f^{\otimes (m-1)}\otimes g^{\otimes (l-1)})\sum_{k=1}^n(Q_{k,n}P_{k,n})\\
&+2(l-1){\mathcal S}_{m+l-2}^n(f^{\otimes m}\otimes g^{\otimes (l-2)})\sum_{k=1}^nP_{k,n}^2\\
&-2^{\frac{3}{2}}m(m-1){\mathcal S}_{m+l-3}^n(f^{\otimes (m-2)}\otimes g^{\otimes (l-1)})\sum_{k=1}^n(Q^2_{k,n}P_{k,n})\\
&-2^{\frac{3}{2}}2m(l-1){\mathcal S}_{m+l-3}^n(f^{\otimes (m-1)}\otimes g^{\otimes (l-2)})\sum_{k=1}^n(Q_{k,n}P^2_{k,n})\\
&-2^{\frac{3}{2}}(l-1)(l-2){\mathcal S}_{m+l-3}^n(f^{\otimes (m-2)}\otimes g^{\otimes (l-1)})\sum_{k=1}^nP^3_{k,n}+\cdots.
\end{split}
\end{equation}

Combining \eqref{S5} and \eqref{S6}, one has
\begin{equation}
\begin{split}
&n^{-\frac{(H+1)(m+l)}{2}}{\mathcal S}_{m+l}^n(f^{\otimes m}\otimes g^{\otimes l})\\
=&n^{-\frac{(H+1)(m+l)}{2}}{\mathcal S}_{m+l-1}^n(f^{\otimes m}\otimes g^{\otimes (l-1)}){\mathcal S}_1^n(g)\\
&-2mn^{-\frac{(H+1)(m+l-2)}{2}}{\mathcal S}_{m+l-2}^n(f^{\otimes(m-1)}\otimes g^{\otimes (l-1)})\left(n^{-(H+1)}\sum_{k=1}^nQ_{k,n}P_{k,n}\right)\\
&-2(l-1)n^{-\frac{(H+1)(m+l-2)}{2}}{\mathcal S}_{m+l-2}^n(f^{\otimes m}\otimes g^{\otimes (l-2)})\left(n^{-(H+1)}\sum_{k=1}^nP^2_{k,n}\right)\\
&+n^{-\frac{(H+1)(m+l)}{2}}R'_{(m+l),n},
\end{split}
\end{equation}
where $R'_{(m+l),n}$ consists of the last 3 lines of \eqref{S6}. It is obvious, as \eqref{2ergodicity}, that the quantities in two parentheses above converge, ${\mathbb Q}$ a.s., to $\frac{1}{2}\langle f,g\rangle_H$ and $\frac{1}{2}\|g\|^2_{H}$, respectively. If we can show $n^{-\frac{(H+1)(m+l)}{2}}R'_{(m+l),n}$ converging weakly to 0 as $n\rightarrow\infty$, then by induction hypothesis \eqref{tensor2}, we obtain 
\begin{equation}\label{S7}
\begin{split}
&n^{-\frac{(H+1)(m+l)}{2}}{\mathcal S}_{m+l}^n(f^{\otimes m}\otimes g^{\otimes l})\longrightarrow I_{m+l-1}(f^{\otimes m}\otimes g^{\otimes (l-1)})I_1(g)\\
-&mI_{m+l-2}(f^{\otimes (m-1)}\otimes g^{\otimes (l-1)})\langle f,g\rangle_H-(l-1)I_{m+l-2}(f^{\otimes m}\otimes g^{\otimes (l-2)})\|g\|^2_{H}.
\end{split}
\end{equation}
By \eqref{contract3}, it follows from \eqref{S7} that
\begin{equation}\label{S8}
\begin{split}
&n^{-\frac{(H+1)(m+l)}{2}}{\mathcal S}_{m+l}^n(f^{\otimes m}\otimes g^{\otimes l})\longrightarrow I_{m+l}(f^{\otimes m}\otimes g^{\otimes l})
\end{split}
\end{equation}
as $n\longrightarrow\infty$, which is what we want.

Now we go back to $n^{-\frac{(H+1)(m+l)}{2}}R'_{(m+l),n}\rightarrow 0$. It suffices to show 
\begin{align*}
2^{\frac{3}{2}}n^{-\frac{(H+1)(m+l)}{2}}{\mathcal S}_{m+l-3}^n(f^{\otimes (m-2)}\otimes g^{\otimes (l-1)})\sum_{k=1}^n(Q^2_{k,n}P_{k,n})\longrightarrow 0.
\end{align*}
since the other terms can be dealt with in a similar manner. Actually, by the law of large number for triangular arrays again as done in \eqref{2ergodicity}, combining \eqref{boundedness}, we have
\begin{align*}
&2^{\frac{3}{2}}n^{-\frac{(H+1)(m+l)}{2}}{\mathcal S}_{m+l-3}^n(f^{\otimes (m-2)}\otimes g^{\otimes (l-1)})\sum_{k=1}^n(Q^2_{k,n}P_{k,n})\\
=&n^{-\frac{(H+1)(m+l-3)}{2}}{\mathcal S}_{m+l-3}^n(f^{\otimes (m-2)}\otimes g^{\otimes (l-1)})2^{\frac{3}{2}}n^{-\frac{3}{2}}\sum_{k=1}^nn^{-\frac{3H}{2}}(Q^2_{k,n}P_{k,n})\longrightarrow 0
\end{align*}
by Slutsky's theorem and induction hypothesis. Finally, we can mimic the above proof to show the result is true for all $f_i's$ are the linear combinations of indicators, thus by Lemma \ref{chaosapproximation} the proof is completed.\qed

\subsection{The convergence of partition function $Z^\om_n(\beta)$.} We go back to the original partition function $Z^\om_n(\beta)$.
Let
\begin{equation}\label{omegan1}
{\om}_n(i,x)=\frac{\e^{\be n^{-\vr}\om(i,x)-\Lambda(\be n^{-\vr})}-1}{\be n^{-\vr}}\stackrel{\De}{=}F(n,\om(i,x)),
\end{equation}
where $\Lambda(\cdot)$ is the Log-Laplace of $\om(i,x)$. Thus, we get a mean zero stationary field ${\om}_n(i,x)$ ($n$-dependent), which is a non-linear functionals of $\om(i,x)$. The covariance of ${\om}_n(i,x)$ and ${\om}_n(i,y)$ is given by 
\begin{align}\label{omegancovariance}
\begin{split}
\EE_\QQ({\om}_n(i,x){\om}_n(i,y))&=\frac{1}{\be^2 n^{-2\vr}}\EE_\QQ\{\e^{\be n^{-\vr}({\om}_n(i,x)+{\om}_n(i,y))-2\Lambda(\be n^{-\vr})}-1\}\\
&=\ga(x-y)(1+o(1)):={\tilde\ga}_n(x-y).
\end{split}
\end{align}
We can expand $F(n,z), z\in\mathbb{R}$, by
$$
F(n,z)=\frac{1}{\be n^{-\vr}}\sum_{k=1}^\infty\frac{(\be n^{-\vr})^k}{k!}A_k(z),
$$
where $A_k(z), k\in\mathbb{N}$, is the system of  Appell polynomials\footnote{The Appell polynomials $\{A_k(z),z\in\RR:k\in \NN\cup\{0\}\}$ for a variable $X$ is defined by
$$
\frac{\e^{sz}}{\EE[\e^{sX}]}=\e^{sz-\Lambda(s)}=\sum_{k=0}^\infty \frac{s^k}{k!}A_k(z),
$$
where $\Lambda(z)$ is the Log-Laplace of $X$. When $X\sim N(0,1)$ its Appell polynomials are just Hermite polynomials. For more details about Appell polynomials see \cite{AVRAM and TAQQU 1987annals}.} related to the distribution of $\om$ with $A_0=1$.  Let $c_k, k\in\NN$, be the expansion coefficients of $F^n$ with respect to Appell system $A_k, k\in\NN$. We remark that the Appell rank, which is the least index $k$ such that $c_k\ne 0$,  of $F{(n,u)}$ is 1 and $c_1=1$. Now by \eqref{omegan1},
we have 
\begin{align}
&\e^{-n\Lambda(\be n^{-\vr})}{\mathbb Z}^\om_n=\e^{-n\Lambda(\be n^{-\vr})}{\EE}_{\PP}\e^{\be n^{-\vr}\sum_{i=1}^n\om(i,S_i)}={\EE}_{\PP}\Pi_{i=1}^n(1+\be n^{-\vr}{\om}_n(i,S_i)).
\end{align}

Consider the new modified partition function $\mathfrak {Z}_n^{{\om}_n}(\be n^{-\vr})$ with ${\om}_n$ replacing the original $\om$: 
\begin{align*}
\mathfrak {Z}_n^{{\om}_n}(\be n^{-\vr})=&{\mathbb E_{\PP}}\left[\prod_{i=1}^n(1+\be n^{-\vr}{\om}_n(i,S_i))\right]\\
=&{\EE_{\PP}}\left[1+\sum_{k=1}^n\be^kn^{-k\vr}\sum_{{\mathbf i}\in D_k^n}\prod_{j=1}^k{\om}_n(i_j,S_{i_j})\right]\\
=&1+\sum_{k=1}^n\be^kn^{-k\vr}\sum_{{\mathbf i}\in D_k^n}\sum_{{\mathbf x}\in {\ZZ}^k}\left[\prod_{j=1}^k{\om}_n(i_j,x_j)p_k({\mathbf i},{\mathbf x})\right],
\end{align*}
and the corresponding weighted $U-$statistics ${\mathcal S}_k^n$ by 
$$
{\mathcal S}_k^n(f,{\om}_n)=2^{k/2}\sum_{{\mathbf i}\in E_k^n}\sum_{{\mathbf x}\in {\mathbb Z}^k}
{\bar f}_n(\frac{{\mathbf i}}{n},\frac{{\mathbf x}}{\sqrt{n}}){\om}_n({\mathbf i},{\mathbf x})\mathbbm{1}_{\{{\mathbf i}\leftrightarrow{\mathbf x}\}}.
$$
In this case, we also have
\begin{thm}\label{omeganconvergence}
Let ${\om}_n$ is defined by \eqref{omegan1}, and assume the existence of Laplacian transformation of $\xi_{ij}$ in a small neighborhood of zero, let $\la(\cdot)$ be the Log-Laplace of $\xi_{ij}$. Then, for $f\in {\mathcal L}_H^{\otimes k}, k\in\NN$,
\begin{align*}
{\mathcal S}_k^n(f,{\om}_n)\stackrel{D}{\longrightarrow} I_k(f), \quad \mbox{as}\quad n\rightarrow \infty.
\end{align*}
\end{thm}
\proof By combining the proof of Theorem \ref{mainthm}, Lemma \ref{chaosapproximation} and the covariance \eqref{omegancovariance}, we only show it holding for $k=1$ and $f$ of the form $f(t,x)={\mathbbm1}_{[t_0,t_1]\times[x_0, x_1]}(t,x)$ for some $0\leq t_0\leq t_1\leq 1, x_0\leq x_1\in\RR$ as before.
 \begin{equation}\label{U1}
 {\mathcal S}_1^n(f,{\om}_n)=2^{1/2}\sum_{{ i}\in E_k^n, nt_0\leq { i}\leq nt_1}\sum_{{x}\in {\mathbb Z},\sqrt{n}x_0\leq {x}\leq\sqrt{n}x_1}
{\om}_n({i},{ x}){\mathbbm 1}_{\{{ i}\leftrightarrow{x}\}}.
\end{equation}

First, we compute the variance of ${\mathcal S}_{1}^n(f,{\om}_n)$. Similarly, we have
\begin{align*}
&{\mathbb E}_{{\mathbb Q}} [({\mathcal S}_{1}^n(f,{\om}_n))^2]\\
=&\la n^{\frac{3}{2}}(t_1-t_0)(x_1-x_0){\ga}_n(0)+2\la n(t_1-t_0)\sum_{k=1}^{N-1}(N-k){\ga}_n(k)+o(n)\\
=&\la n^{\frac{3}{2}}(t_1-t_0)(x_1-x_0)({\ga}(0)(1+o(1)))+2\la n(t_1-t_0)\sum_{k=1}^{N-1}(N-k)(\ga(k)(1+o(1)))+o(n).
\end{align*}
It is obvious that 
\begin{align*}
&n^{-(H+1)}{\mathbb E}_{{\mathbb Q}} [({\mathcal S}_{1}^n(f))^2]\longrightarrow \frac{\la (t_1-t_0)(x_1-x_0)^{2H}}{H(2H-1)}={\mathbb E}_H\left(\int_0^1\int_Rf(t,x)W(\D t\D x)\right)^2.
\end{align*}

Next we show the asymptotic normal  of ${\mathcal S}_{1}^n(f,{\om}_n)$. The quantity ${\om}_n$ in the right hand side of \eqref{U1} is a nonlinear functionals of the moving process $\omega$.  For the theories of central or non-central limit theorem of stationary process one can be referred to \cite{Dobrushin1979Non} for Gaussian case and \cite{AVRAM and TAQQU 1987annals, Giraitis1985Central,Surgailis1982Zones} for non-Gaussian case, which the Apell expansions for non-linear functionals are used extensively. Here we will follow the procedure used in \cite{Ho and Hsing1998Annals}, where a different approach without polynomials expansion was adopted. Here we  verify that the functional ${\om}_n$ satisfies the conditions of Corollary 3.3 in \cite{Ho and Hsing1998Annals}. 

For this sake, we introduce some similar notations as in \cite{Ho and Hsing1998Annals}. Define $\om ( i , x)(0) = 0$  and 
$$
\om ( i , x)(y) = \sum _ { |z| < |y| } \psi _ { z -x} \xi _ {i, z } , \quad \tilde { \om }( i , x)(y) = { \om }( i , x)- \om ( i , x)(y)  , \quad y \in\NN.
$$
Let 
$$
F_{i,x,y}(n,u)=\EE[F(n,u+\om ( i , x)(y))]=\frac{1}{a}[\exp\{au-\sum_{|z|\geq y}\lambda(a\psi_{x-z})\}-1],
$$
where $a=\beta n^{-\rho}$. For $l=0,1,2,3$, the $l$th derivative  $F^{(l)}_{i,x,y}(n,u)$ of $F_{i,x,y}(n,u)$ exists and for all $V>0$
$$
F^{(l)}_{i,x,y,V}(n,u)=\sup_{|v|\leq V}|F^{(l)}_{i,x,y}(n,u+v)|,
$$
has an explicit expression, which is continuous in $u$, e.g., 
$$
F^{(3)}_{i,x,y,V}(n,u)=a^2\exp\{au+aV-\sum_{|z|\geq y}\lambda(a\psi_{x-z})\}.
$$
Furthermore, we have, for all $u\in\RR$,
$$
\sup_{I\subset \ZZ}\EE[F^{(l)}_{i,x,y,V}(n,u+\sum_{y\in I}\psi_y\xi_{i,y})]^4<\infty
$$
by the assumption on the exponential integrability of $\xi_{i,y}$ in a small neighbor of zero since $n$ is large enough, where the sup is taken over all subsets $I$ of $\ZZ$. These are called $C(l,y,V)$ conditions. Since the Appell rank of $F{(n,u)}$ is 1 and $c_1=1$ we can decompose the rhs of \eqref{U1} into two parts as follows:
$$
\begin{aligned}
&2^{1/2}\sum_{{ i}\in E_k^n, nt_0\leq { i}\leq nt_1}\sum_{{x}\in {\mathbb Z},\sqrt{n}x_0\leq {x}\leq\sqrt{n}x_1}
[{\om}_n({i},{ x})-\om(i,x)]1_{\{{ i}\leftrightarrow{x}\}}\\
&+2^{1/2}\sum_{{ i}\in E_k^n, nt_0\leq { i}\leq nt_1}\sum_{{x}\in {\mathbb Z},\sqrt{n}x_0\leq {x}\leq\sqrt{n}x_1}
\om(i,x)1_{\{{ i}\leftrightarrow{x}\}}=\mathcal{R}_n+\mathcal{M}_n.
\end{aligned}
$$
We have shown that 
$$
n^{-\frac{H+1}{2}}\mathcal{M}_n\stackrel{D}{\longrightarrow}\int_0^1\int_Rf(t,x)W(\D t\D x)
$$
as $n\to\infty$. According to Theorem 3.1 or Corollary 3.3 in \cite{Ho and Hsing1998Annals}, we know 
$$
[\text{Var}_\QQ({\mathcal S}_{1}^n(f,{\om}_n))]^{-1}\EE_\QQ(\mathcal{R}_n)^2\longrightarrow 0
$$
as $n\to \infty$. Hence, by Slutsky's theorem, we have the desired result of this theorem.\qed
 
 Finally, we arrive at the following result of this paper.
 \begin{thm}
 \label{convergence12}
Let $\{u(t,x),(t,x)\in [0,1]\times{\mathbb R}\}$ be the solution to \eqref{SHE1} with parameter $\sqrt{2}\beta$, initial data $u(x)=\de(x)$. And let $Z^\om_n$ be the partition function \eqref{partition1} of random polymer in the random environment $\{\om(n,x):n\geq 0, x\in\ZZ\}$ with the representation \eqref{envir1}. Then
\begin{align*}
&\frac{\sqrt{n}}{2}\e^{-nt\Lambda(\be n^{-\vr})}Z^\om_n(nt,\sqrt{n}x;\be n^{-\vr})\longrightarrow u(t,x) \quad \mbox{weakly, as}\quad n\longrightarrow\infty.
\end{align*}
\end{thm}
 
\section{Tightness}

In this section we prove the approximation process 
$$z_n(s,y,t,x):=\sqrt{n}\mathfrak {Z}^\om_n(ns,\sqrt{n}y;nt,\sqrt{n}x;\be n^{-\vr})
$$
 is tight by Kolmogorov's criterion. Here we also consider only the case of two-parameter, i.e., tightness of $z_n(t,x)=\sqrt{n}\mathfrak {Z}^\om_n(nt,\sqrt{n}x;\be n^{-\vr})$. For four-parameter field $z_n(s,y,t,x)$ see the Remark in \cite[Section 5]{alberts2014}. Also $z_n(t,x)$ is designed to be jointed pitch by pitch. To be concrete, on every rectangle $\big(\frac{{i}-1}{n},\frac{{ i}}{n}\big]\times\big(\frac{{x}-1}{\sqrt{n}},\frac{{x}+1}{\sqrt{n}}\big]$ for $i\in\NN, x\in\ZZ$, the roof of $z_n(t,x)$ is flat or pasted by two triangles. 

Since the components of environment in the time direction are independent, the Markovian property still holds. Now we mimic the procedure as in \cite{alberts2014} to obtain the difference equation for $\sqrt{n}\mathfrak {Z}^\om_n$. In accordance with the definition of $\mathfrak {Z}^\om_n$, one has, by condition on time $k$,
\begin{align*}
&{\mathfrak {Z}}^\om_n(k+1,x;\beta)-{\mathfrak {Z}}^\om_n(k,x;\beta)\\
=&\frac{1}{2}\De {\mathfrak Z}^\om_n(k,x;\beta)+\be\om(k+1,x){\bar {\mathfrak Z}}^\om_n(k,x;\be),
\end{align*}
which is a discrete version of heat equation, where ${\bar {\mathfrak Z}}^\om_n(k,x;\be)=\frac{1}{2}[{\mathfrak Z}^\om_n(k+1,x;\beta)+{\mathfrak Z}^\om_n(k-1,x;\beta)]$ and $\De$ is the discrete Laplacian operator. 
Therefore, by Duhamel's principle, it results in the following equation:
\begin{align*}
&{\mathfrak {Z}}^\om_n(k,x;\beta)=p(k,x)+\beta\sum_{i=1}^k\sum_{y}\om(i,y){\bar {\mathfrak Z}}^\om_n(i-1,y;\be)p(k-i,x-y).
\end{align*}
Now, by scaling ${{\mathfrak Z}}^\om_n(k,x;\be)$, we have
\begin{align}\label{eqzn1}
z_n(t,x)=p_n(t,x)+n^{-\frac{3}{2}}\be\sum_{\substack{s\in[0,t]\cap n^{-1}{\mathbb Z}\\ y\in n^{-1/2}{\mathbb Z}}} p_n(t-s,x-y)\bar{z}_n(s,y)\om_n(s,y),
\end{align} 
where $\om_n(s,y)=n^{1-\varrho}\om(ns,\sqrt{n}y)$, $p_n(t,x)=\sqrt{n}p(\lfloor nt\rfloor, \sqrt{n}x)$. 

For large enough $n\in{\mathbb N}$, as in \cite{Dalang1999}, we define random martingale measure $M_n$ by 
\begin{align*}
M_n([0,t]\times A)=\sum_{\substack{s\in[0,t]\cap n^{-1}{\mathbb Z}\\ y\in n^{-1/2}{\mathbb Z}\cap A}} \om_n(s,y)\stackrel{\De}{=}M_n(t,A)
\end{align*} 
 with $A\in {\mathscr A}=\{A\in {\mathcal B}(\RR): |A|_H<\infty\}$, where $|A|_H=\int_A\int_AK(u,v)\D u\D v$. Then for $A\in {\mathscr A}$ fixed, $M_n(\cdot,A)=\{M_n(t, A), t\in n^{-1}\ZZ\}$ is a stationary independent increments process. Let ${\mathscr P}$ be the predictable $\si$-algebra generated by $M_n,n=1,2,\dots$. By the definition of ${\bar {\mathfrak Z}}^\om_n(i-1,y;\be)$, we know that $\bar{z}_n(s,\cdot)$ is predictable process. Hence the sum 
$$
n^{-\frac{3}{2}}\sum_{\substack{s\in[0,t]\cap n^{-1}{\mathbb Z}}}\sum_{ y\in n^{-1/2}{\mathbb Z}}\bar{z}_n(s,y)\om_n(s,y)\stackrel{\De}{=}\int_0^t\int_{\RR}\bar{z}_n(s,y)\om_n(s,y)\D s\D y
$$ 
is a martingale with $\si$-algebra ${\mathcal F}_t, t\geq 0$, generated by $\om$ up to time $t$.
 Its quadratic variation process is given by 
 $$
 n^{-1-2\vr}\sum_{\substack{s\in[0,t]\cap n^{-1}{\mathbb Z}}}\sum_{ y_1\in n^{-1/2}{\mathbb Z}}\sum_{ y_2\in n^{-1/2}{\mathbb Z}}\bar{z}_n(s,y_1)\ga(\sqrt{n}y_1-\sqrt{n}y_2)\bar{z}_n(s,y_2).
 $$
Then the second term of the r.h.s. of \eqref{eqzn1} can be understood in the sense of the stochastic integral of the kernel of simple random walk with respect to the martingale
 $$\int_0^t\int_{\RR}\bar{z}_n(s,y)\om_n(s,y)\D s\D y.$$ 
 Therefore we can rewrite \eqref{eqzn1} as the integral form
\begin{align}\label{eqzn2}
z_n(t,x)=p_n(t,x)+\be\int_0^t\int_\RR p_n(t-s,x-y)z_n(s,y)\om_n(s,y)\D s\D y.
\end{align} 
 
  We now turn back to 
\eqref{eqzn1} or \eqref{eqzn2} and check the Kolmogorov's criterion for $z_n$. To this end, let $q\geq 1, \iota>0,\tau>0$, which will be specified later. In fact, we know for any $q>0$, there exists a  constant $C_q$ such that for any 
\begin{equation}
\max_{t\in [0,1],x\in\RR}\EE[|z_n(t,x)|^{2q}]\leq C_q,
\end{equation}
which can be proved by routine method adopted in SPDE, such as Burkholder inequality, Young inequality (See the estimation for $Q_2$ below). In what follows, $C$ or $C_q$ is a generic constant independent of $n,s,x,$ etc, and may be changed from line to line.

Now for $t_1>t_2\in[0,1], x\in\RR$, we have
\begin{align}\label{eq time}
\begin{split}
&{\mathbb E}|z_n(t_1,x)-z_n(t_2,x)|^{2q}\\
\leq &C_q{\mathbb E}\left|\int_{t_2}^{t_1}\int_{\mathbb R}p_n(t_1-s,x-y)z_n(s,y)\om_n(s,y)\D s\D y\right|^{2q}\\
&+C_q{\mathbb E}\left|\int_0^{t_2}\int_{\mathbb R}(p_n(t_1-s,x-y)-p_n(t_2-s,x-y))z_n(s,y)\om_n(s,y)\D s\D y\right|^{2q}\\
:=&C_qQ_1+C_qQ_2.
\end{split}
\end{align}

By Burkholder inequalities, we have that $Q_2$ is less than, up to a constant multiplier, 
\begin{align}\label{Q2}
\begin{split}
&n^{-Hq-q}{\mathbb E}\Big(\!\sum_{\{s\in[0,t_2]\cap \frac{1}{n}{\mathbb Z}\}}\sum_{\{ y_1,y_2\in \frac{1}{\sqrt{n}}{\mathbb Z}\}}(p_n(t_1-s,x-y_1)-p_n(t_2-s,x-y_1))\bar{z}_n(s,y_1)\\
&\times\ga(\sqrt{n}y_1-\sqrt{n}y_2)\bar{z}_n(s,y_2)(p_n(t_1-s,x-y_2)-p_n(t_2-s,x-y_2))\Big)^{q}.
\end{split}
\end{align}
The expectation above of $q$-power of the sum can be written as
\begin{align*}
\sum_{\{s_1\in[0,t_2]\cap \frac{1}{n}{\mathbb Z}\}}&\sum_{\{ y_{1,1},y_{2,1}\in \frac{1}{\sqrt{n}}{\mathbb Z}\}}\cdots\sum_{\{s_q\in[0,t_2]\cap \frac{1}{n}{\mathbb Z}\}}\sum_{\{ y_{1,i},y_{2,i}\in \frac{1}{\sqrt{n}}{\mathbb Z}\}}\EE(\Pi_{i=1}^q\bar{z}_n(s_i,y_{1,i})\bar{z}_n(s_i,y_{2,i}))\\
\times&\Pi_{i=1}^q[(p_n(t_1-s_i,x-y_{1,i})-p_n(t_2-s,x-y_{1,i}))\\
&~~~~~~~~~~\times\ga(\sqrt{n}y_{1,i}-\sqrt{n}y_{2,i})(p_n(t_1-s_i,x-y_{2,i})-p_n(t_2-s_i,x-y_{2,i}))]
\end{align*}
By the generalized H\"{o}lder inequality,
\begin{align*}
\EE(\Pi_{i=1}^q\bar{z}_n(s_i,y_{1,i})\bar{z}_n(s_i,y_{2,i}))&\leq \Pi_{i=1}^q(\EE|\bar{z}_n(s_i,y_{1,i})|^{2q})^{1/2q}(\EE|\bar{z}_n(s_i,y_{2,i})|^{2q})^{1/2q}\\
&\leq C\Pi_{i=1}^q(\EE|{z}_n(s_i,y_{1,i})|^{2q})^{1/2q}(\EE|{z}_n(s_i,y_{2,i})|^{2q})^{1/2q}.
\end{align*}
Consequently, \eqref{Q2} is bounded above by
\begin{align*}
Cn^{-Hq-q}\sum_{\{s_1\in[0,t_2]\cap \frac{1}{n}{\mathbb Z}\}}&\sum_{\{ y_{1,1},y_{2,1}\in \frac{1}{\sqrt{n}}{\mathbb Z}\}}\cdots\sum_{\{s_q\in[0,t_2]\cap \frac{1}{n}{\mathbb Z}\}}\sum_{\{ y_{1,q},y_{2,q}\in \frac{1}{\sqrt{n}}{\mathbb Z}\}}\\
&\Pi_{i=1}^q(\EE|{z}_n(s_i,y_{1,i})|^{2q})^{1/2q}(\EE|{z}_n(s_i,y_{2,i})|^{2q})^{1/2q}\\
\times&\Pi_{i=1}^q(p_n(t_1-s_i,x-y_{1,i})-p_n(t_2-s,x-y_{1,i}))\\
&~~~~\times\ga(\sqrt{n}y_1-\sqrt{n}y_2)(p_n(t_1-s_i,x-y_{2,i})-p_n(t_2-s_i,x-y_{2,i}))
\end{align*}
By the positivity, the quantity above is less than
\begin{align*}
&C\max_{s,x}{\mathbb E}(z_n^{2q}(s,y))\bigg(n^{-H-1}\sum_{\{s\in[0,t_2]\cap \frac{1}{n}{\mathbb Z}\}}\sum_{\{ y_1,y_2\in \frac{1}{\sqrt{n}}{\mathbb Z}\}}\\
&~~~~~~~~~~~~~~~(p_n(t_1-s,x-y_1)-p_n(t_2-s,x-y_1))\\
&~~~~~~~~~~\times\ga(\sqrt{n}y_1-\sqrt{n}y_2)(p_n(t_1-s,x-y_2)-p_n(t_2-s,x-y_2))\bigg)^{q}
\end{align*}
\begin{align*}
=&C\max_{s,x}{\mathbb E}(|z_n(s,y)|^{2q})\bigg(n^{-H-1}\int_{-\pi}^\pi \!\!G(\D\eta)\!\!\!\sum_{\{s\in[0,t_1]\cap \frac{1}{n}{\mathbb Z}\}}\sum_{\{ y_1\in \frac{1}{\sqrt{n}}{\mathbb Z}\}}\\
&~~~~~~~~~~~~~~~(p_n(t_1-s,x-y_1)-p_n(t_2-s,x-y_1))\e^{\i\sqrt{n}y_1\eta}\\
&~~~~~~~~~~\times\sum_{\{ y_2\in \frac{1}{\sqrt{n}}{\mathbb Z}\}}(p_n(t_1-s,x-y_2)-p_n(t_2-s,x-y_2))\e^{-\i\sqrt{n}y_2\eta}\bigg)^{q},
\end{align*}
where $G(\D\eta)$ is the spectrum measure of $\ga$. Noticing that
\begin{align*}
&\sum_{\{ y_1\in \frac{1}{\sqrt{n}}{\mathbb Z}\}}(p_n(t_1-s,x-y_1)-p_n(t_2-s,x-y_1))\e^{\i\sqrt{n}y_1\eta}\\
=&\sqrt{n}\e^{\i\sqrt{n}x\eta}((\frac{1}{2}\e^{\i\eta}+\frac{1}{2}\e^{-i\eta})^{n(t_1-s)}-(\frac{1}{2}\e^{\i\eta}+\frac{1}{2}\e^{-\i\eta})^{n(t_2-s)}),
\end{align*}
we have that $Q_2$ is controlled by
\begin{align*}
&\max_{s,x}{\mathbb E}(|z_n(s,y)|^{2q})\bigg|n^{-H}\int_{-\pi}^\pi \!\!G(\D\eta)\!\!\!\sum_{\{s\in[0,t_2]\cap \frac{1}{n}{\mathbb Z}\}}\\
&~~~~~~~~~~~~~~~~~~~~~~\Big((\frac{1}{2}\e^{\i\eta}+\frac{1}{2}\e^{-\i\eta})^{n(t_1-s)}-(\frac{1}{2}\e^{\i\eta}+\frac{1}{2}\e^{-\i\eta})^{n(t_2-s)}\Big)^2\bigg|^{q}\\
\leq&\max_{s,x}{\mathbb E}(|z_n(s,y)|^{2q})\bigg|n^{-H}\int_{-\pi}^\pi \!\!G(\D\eta)\!\!\!\sum_{\{s\in[0,t_2]\cap \frac{1}{n}{\mathbb Z}\}}(\cos\eta)^{2n(t_2-s)}\Big((\cos\eta)^{n(t_1-t_2)}-1\Big)^2\bigg|^{q}\\
\leq&\max_{s,x}{\mathbb E}(|z_n(s,y)|^{2q})\bigg|n^{-1}\int_{-\sqrt{n}\pi}^{\sqrt{n}\pi} \!\!n^{1-H}G(\frac{\D\eta}{\sqrt{n}})\frac{1-(\cos\frac{\eta}{\sqrt{n}})^{2nt_2+2}}{\sin^2\frac{\eta}{\sqrt{n}}}\Big((\cos\frac{\eta}{\sqrt{n}})^{n(t_1-t_2)}-1\Big)^2\bigg|^{q}.
\end{align*}
Changing variables by $n(t_1-t_2)=k$, then $\sqrt{t_1-t_2}\eta=z$, we have
\begin{align*}
Q_2
\leq&\max_{s,x}{\mathbb E}(|z_n(s,y)|^{2q})\bigg|(t_1-t_2)^Hk^{-H}\int_{-\sqrt{k}\pi}^{\sqrt{k}\pi} \!\!G(\frac{\D z}{\sqrt{k}})\frac{1-(\cos\frac{z}{\sqrt{k}})^{\frac{2kt_2}{t_1-t_2}+2}}{\sin^2\frac{z}{\sqrt{k}}}\Big((\cos\frac{z}{\sqrt{k}})^{k}-1\Big)^2\bigg|^{q}.
\end{align*}

By the similar manner as \eqref{E1}, it follows that 
\begin{align*}
Q_2
\leq&\max_{s,x}{\mathbb E}(|z_n(s,y)|^{2q})\bigg(\int_{-\infty}^{\infty} \frac{\big(1-\e^{-\frac{(t_1-t_2)z^2}{t_2}}\big)\big(\e^{-\frac{z^2}{2}}-1\big)^2}{|z|^{1+H}}\D z\bigg)^{q}(t_1-t_2)^{qH} \leq C(t_1-t_2)^{qH},
\end{align*}
since the integral exists obviously.

In a similar way, we have
\begin{align*}
Q_1
\leq&\max_{s,x}{\mathbb E}(|z_n(s,y)|^{2q})\bigg(\int_{-\infty}^{\infty} \frac{\big(1-\e^{-\frac{(t_1-t_2)z^2}{t_2}}\big)\big(\e^{-\frac{z^2}{2}}-1\big)^2}{|z|^{1+H}}\D z\bigg)^{q}(t_1-t_2)^{qH} \leq C(t_1-t_2)^{qH}.
\end{align*}

Now we consider the moment of spatial increment ${\mathbb E}|z_n(t,x)-z_n(t,y)|^{2q}$. Actually, two terms we are concerned. One is $Q_3:=|p_n(t,x)-p_n(t,y)|^{2q}\leq C|x-y|^{\iota q}$, for some $0<\iota<H$ and $C>0$, which can be estimated easily by local central limit theorem (see \cite[Proposition 2.4.1]{Lawler2010Random}) under condition of $t\in[\epsilon, 1]$, $\epsilon>0$ fixed. The second one is, as \eqref{Q2},

\begin{align}\label{Q4}
\begin{split}
&n^{-Hq-q}{\mathbb E}\Big(\!\sum_{\{s\in[0,t]\cap \frac{1}{n}{\mathbb Z}\}}\sum_{\{ y_1,y_2\in \frac{1}{\sqrt{n}}{\mathbb Z}\}}(p_n(t-s,x-y_1)-p_n(t-s,y-y_1))\bar{z}_n(s,y_1)\\
&\times\ga(\sqrt{n}y_1-\sqrt{n}y_2)\bar{z}_n(s,y_2)(p_n(t-s,x-y_2)-p_n(t-s,y-y_2))\Big)^{q}:=Q_4.
\end{split}
\end{align}
Similarly, one can show 
\begin{align*}
Q_4&\leq C\max_{s,x}{\mathbb E}(|z_n(s,y)|^{2q})n^{-qH}\Big\{\int_{-\pi}^\pi(\e^{\imath\sqrt{n}x\eta}-\e^{\imath\sqrt{n}y\eta})^2\frac{1-(\cos^2\eta)^{nt+1}}{\sin^2\eta}G(\D\eta)\Big\}^q\\
&\leq C\max_{s,x}{\mathbb E}(|z_n(s,y)|^{2q})n^{-q}\Big\{\int_{-\sqrt{n}\pi}^{\sqrt{n}\pi}(\e^{\imath x\eta}-\e^{\imath y\eta})^2\frac{1-(\cos^2\frac{\eta}{\sqrt{n}})^{nt+1}}{\sin^2\frac{\eta}{\sqrt{n}}}n^{1-H}G(\frac{\D\eta}{\sqrt{n}})\Big\}^q.
\end{align*}
By inequality $|\sin(x-y)\eta|\leq|x-y|^r|\eta|^r$, for some $0<r<H$, we have
\begin{align*}
Q_4&\leq C\max_{s,x}{\mathbb E}(|z_n(s,y)|^{2q})|x-y|^{2rq}\Big\{\int_{-\infty}^{\infty}\frac{1-\e^{-t\eta^2}}{\eta^{1+2H-2r}}\D \eta\Big\}^q.
\end{align*}

In conclusion, combining $Q_1,Q_2,Q_3$ and $Q_4$, we finally get

\begin{thm}\label{tight1}
Let $\epsilon>0$ be small enough. For any $n\in\NN$, $t,s\in[\epsilon,1]$ and $x,y\in\RR$, for some $q>1$, there exist constant $C_\epsilon>0,0<\iota<H$, such that
\begin{align}\label{tight}
\EE|z_n(t,x)-z_n(s,y)|^{2q}\leq C_\epsilon(|t-s|^{Hq}+|x-y|^{\iota q}).
\end{align}
Moreover, if $2q$-order moment of $\om$ is finite for $q>\frac{2}{H}$, then the family of process $\{z_n\}_{n=1}^\infty$ is tight in $C([\epsilon,1],\RR)$.
\end{thm}
\proof \eqref{tight} is the consequence of the above computation. The tightness is obtained by applying Kolmogorov's tightness criterion if we chose $q$ so that $q\iota>1$. The existence of such $q$ is guaranteed by the assumption that $2q$ moment of $\om$ is finite.\qed 

\section*{ Acknowlegement}  This work  is supported by the National Natural Science Foundation of China (Grants no. 11571262, 11731012 and 11971361). The author thanks Fuqing Gao for much stimulating discussion and thanks Yingxia Chen for her thorough reading the draft and spotting some mistakes. Also the author would appreciate the referees whose comments help improve the paper greatly. 

\end{document}